\documentclass[10pt]{article}
\usepackage{amscd,amsmath,amssymb,amsthm,mathrsfs,mathptm}
\usepackage[all]{xy}
\usepackage{color}
\usepackage{multirow}

\objectmargin{1pt}

\newtheorem{Def}{Definition}[section]
\newtheorem{Prop}[Def]{Proposition}
\newtheorem{Theo}[Def]{Theorem}
\newtheorem{Lem}[Def]{Lemma}
\newtheorem{Koro}[Def]{Corollary}

\newcommand{\E}[4]{{\rm E}_{#1}^{#2}(#3, #4)}
\newcommand{\Ex}[3]{{\rm E}_{#1}^{#2}(#3)}

\newcommand{\add}{{\rm add}}

\newcommand{\Hom}{{\rm Hom }}
\newcommand{\rad}{{\rm rad}}

\newcommand{\soc}{{\rm soc}}
\renewcommand{\top}{{\rm top}}
\newcommand{\StHom}{{\rm \underline{Hom} }}

\newcommand{\C}[1]{\mathscr{C}(#1)}

\newcommand{\Cb}[1]{{\mathscr{C}^b}(#1)}

\newcommand{\Ker}{{\rm Ker}}
\newcommand{\cpx}[1]{#1^{\bullet}}
\newcommand{\HomP}{{\rm Hom}^{\bullet}}
\newcommand{\D}[1]{\mathscr{D}(#1)}

\newcommand{\Db}[1]{ \mathscr{D}^{\rm b}(#1)}
\newcommand{\End}{{\rm End}}
\newcommand{\Ext}{{\rm Ext}}
\newcommand{\K}[1]{\mathscr{K}(#1)}

\newcommand{\Kb}[1]{ \mathscr{K}^{\rm b}(#1)}
\newcommand{\modcat}[1]{#1\mbox{{\rm -mod}}}
\newcommand{\Modcat}[1]{#1\mbox{{\rm -Mod}}}
\newcommand{\stmodcat}[1]{#1\mbox{{\rm -{\underline{mod}}}}}
\newcommand{\pmodcat}[1]{#1\mbox{{\rm -proj}}}
\newcommand{\imodcat}[1]{#1\mbox{{\rm -inj}}}
\newcommand{\opp}{^{\rm op}}
\newcommand{\otimesL}{\otimes^{\rm\bf L}}
\newcommand{\otimesP}{\otimes^{\bullet}}

\newcommand{\stp}[1]{\nu_{#1}\mbox{\rm -stp}}

\renewcommand{\leq}{\leqslant}
\renewcommand{\geq}{\geqslant}
\newcommand{\lra}{\longrightarrow}
\newcommand{\ra}{\rightarrow}
\newcommand{\lraf}[1]{\stackrel{#1}{\lra}}
\newcommand{\raf}[1]{\stackrel{#1}{\ra}}

\begin{document}

\title{Derived equivalences and stable equivalences of Morita type, II}
\author{ WEI HU {\small AND}  CHANGCHANG XI$^*$}
\date{}
\maketitle

\begin{abstract}
Motivated by understanding the Brou\'e's abelian defect group conjecture from algebraic point of view, we consider the question of how to lift a stable
equivalence of Morita type between arbitrary finite dimensional algebras to a derived equivalence. In this paper, we present
a machinery to solve this question for a class of stable equivalences of Morita type. In particular, we show that every
stable equivalence of Morita type between Frobenius-finite algebras
over an algebraically closed field can be lifted to a derived
equivalence. %Thus Frobenius-finite algebras share many common
%invariants of both derived equivalences and stable equivalences.
Especially, Auslander-Reiten conjecrure is true for stable equivalences of Morita type between Frobenius-finite algebras without semisimple direct summands. Examples of such a class of algebras are abundant, including Auslander algebras, cluster-tilted algebras and certain Frobenius extensions. As a byproduct of our methods, we further show
that, for a Nakayama-stable idempotent element $e$ in an algebra $A$
over an arbitrary field, each tilting complex over $eAe$ can be
extended to a tilting complex over $A$ that induces an almost
$\nu$-stable derived equivalence studied in the first paper of this
series. Moreover, we demonstrate that our techniques are applicable to verify the Brou\'{e}'s abelian defect group conjecture for several cases mentioned by Okuyama.
\end{abstract}

\renewcommand{\thefootnote}{\alph{footnote}}
\setcounter{footnote}{-1} \footnote{ $^*$ Corresponding author.
Email: xicc@cnu.edu.cn; Fax: 0086 10 68903637.}
\renewcommand{\thefootnote}{\alph{footnote}}
\setcounter{footnote}{-1} 
\renewcommand{\thefootnote}{\alph{footnote}}
\setcounter{footnote}{-1} 

\tableofcontents

\section{Introduction}
Derived and stable equivalences of algebras (or categories) are two kinds of fundamental equivalences
both in the representation theory
of algebras and groups  and in the theory of triangulated categories. They preserve many
significant algebraic, geometric or numeric invariants, and provide
surprising and useful new applications to as well as connections with other fields
(see \cite{RickDstable}, \cite{RickDFun} and \cite{Broue1994}). But what are the interrelations between these two classes of equivalences? Rickard showed in \cite{RickDstable} that, for self-injective algebras, derived equivalences imply stable equivalences of Morita type. Conversely, Asashiba proved in \cite{AsashibaLift} that, for representation-finite self-injective algebras, almost every stable equivalence lifts to a derived equivalence. For general algebras, however, little is known about their relationship. That is, one does not know any methods with which one such  equivalence can be constructed from the other
for arbitrary algebras. In \cite{HuXi3}, we started with discussing this kind of questions and
gave methods to construct stable equivalences of Morita type from
almost $\nu$-stable derived equivalences, which generalizes the above-mentioned result of Rickard. Of particular interest is also the converse question: How to get derived equivalences from stable equivalences
of Morita type? A motivation behind this question is the \emph{Brou\'e's Abelian Defect
Group Conjecture} which says that the module categories of a block of
a finite group algebra and its Brauer correspondent should have
equivalent derived categories if their defect groups are abelian. Note that block algebras are self-injective. So the result of Rickard implies that Brou\'e's
conjecture would predicate actually  a stable
equivalence of Morita type between the two block algebras. Also, one knows that
stable equivalences of Morita type between block algebras occur very often in Green correspondences. Another motivation is the
\emph{Auslander-Reiten Conjecture} on stable equivalences which states that two stably equivalent algebras should have the same numbers of non-isomorphic non-projective simple modules (see, for instance, \cite[Conjecture (5), p.409]{AusReiten}). This conjecture is even open for stable equivalences of Morita type.
However, it is valid for those stable equivalences of Morita type that can be lifted to derived equivalences since derived
equivalences preserve the numbers of simple modules, while stable equivalences of Morita type between algebras without semisimple summands preserve the numbers of projective simples.
%For stable equivalences of Morita type, this conjecture is related to the Hochschild homology in degree zero (\cite{lzz}),
%while derived equivalences preserve all Hochschild homology groups.
Thus, the above question is of great interest and we restate
it in purely algebraic point of view.

\medskip
{\bf Question.} {\it Given a stable equivalence of Morita type between
arbitrary finite-dimensional algebras $A$ and $B$ over a field, under which conditions can we construct a derived equivalence
therefrom between $A$ and $B$?}

\medskip
%This question, unfortunately, is still to be understood. In particular, an algebraic approach to the question is missing indeed.

In this paper, we shall provide several answers to this question. Our method developed here is different from the one in \cite{AsashibaLift, HuXi2} and can be used to verify the Brou\'e's Abelian Defect Group Conjecture in some cases (see the last section of the paper).

Our first main result provides a class of algebras, called Frobenius-finite algebras, for which every stable equivalences of Morita type induces a derived equivalence (see Subsection \ref{subsect2.2} for definitions).
Roughly speaking, a Frobenius part of a finite-dimensional algebra $A$ is the largest algebras of the form $eAe$ with $e$ an idempotent element such that $\add(Ae)$ is stable under the Nakayama functor.
An algebra is said to be \emph{Frobenius-finite} if its Frobenius part is a representation-finite algebra. Examples of Frobenius-finite algebras are abundant and capture many interesting class of algebras: Representation-finite algebras, Auslander algebras and cluster-tilted algebras. Also, they can be constructed from triangular matrix rings and Frobenius extensions (for more details and examples see Section \ref{sect5.1}).

\begin{Theo}
Let $k$ be an algebraically closed field. Suppose that $A$ and $B$ are two
finite-dimensional $k$-algebras without semisimple direct summands.
If $A$ is Frobenius-finite, then every individual stable equivalence
of Morita type between $A$ and $B$ lifts to an iterated almost
$\nu$-stable derived equivalence. \label{ThmRepFin}
\end{Theo}

Thus the
class of Frobenius-finite algebras shares many common algebraical and
numerical invariants of derived and stable equivalences. Moreover, Theorem \ref{ThmRepFin} not only extends a result of Asashiba in \cite{AsashibaLift} (in a different direction) to a great context, namely every stable equivalence of Morita type between \emph{arbitrary} representation-finite (not necessarily self-injective) algebras lifts to a derived equivalence, but also provides a method to construct a class of derived equivalences between algebras and their subalgebras because under some mild conditions each stable equivalence of Morita type can be realised as an Frobenius extension of algebras by a result in \cite[Corollary 5.1]{DugasVilla}.

As an immediate consequence of Theorem \ref{ThmRepFin}, we get the following: Let $A$ and $B$ be Frobenius-finite
$k$-algebras over an algebraically closed field and without semisimple direct summands. If they are stably equivalent of Morita type, then $A$ and $B$ have the same number of non-isomorphic, non-projective simple modules.
Here, we do not assume that both $A$ and $B$ have no nodes, comparing with a result in \cite{mv}.

Recall that a finite-dimensional $k$-algebra $A$ is called an
\emph{Auslander algebra} if it has global dimension at most $2$ and
dominant dimension at least $2$. Algebras of global dimension at most $2$ seem of great interest in representation theory because they are quasi-hereditary (see \cite{dr}) and every finite-dimensional algebra (up to Morita equivalence) can be obtained from an algebra of global dimension $2$ by universal localization
(see \cite{nrs}).
Since Auslander algebras and cluster-tilted algebras are Frobenius-finite, we have the following consequence of Theorem \ref{ThmRepFin}.

\begin{Koro} Suppose that $A$ and $B$ are  finite-dimensional algebras over algebraically closed field and without semisimple direct summands. If $A$ is an Auslander algebra or a cluster-tilted algebra, then every individual stable equivalence of Morita type between $A$ and $B$ lifts to an iterated almost $\nu$-stable derived equivalence.
\end{Koro}

Our second main result which lays a base for the  proof of Theorem \ref{ThmRepFin} provides a general criterion for lifting a stable equivalence of Morita type to an iterated almost $\nu$-stable derived equivalence. Though this criterion looks technical,  it is more suitable for applications to
surgery. Recall that an idempotent $e$ of an algebra $A$ is said to be \emph{$\nu$-stable }if $\add(\nu_AAe)= \add(Ae)$, where $\nu_A$ is the Nakayama functor of $A$.

\begin{Theo}
Let $A$ and $B$ be finite-dimensional algebras over a field and without semisimple
direct summands, such that $A/\rad(A)$ and $B/\rad(B)$ are separable.
Let $e$ and $f$ be two $\nu$-stable idempotent elements in $A$ and $B$,
respectively, and let $\Phi: \stmodcat{A}\ra\stmodcat{B}$ be a stable
equivalence of Morita type between $A$ and $B$. Suppose that $\Phi$
satisfies the following two conditions:

\smallskip
$(1)$ For all simple $A$-modules $S$ with $e\cdot S=0$, $\Phi(S)$ is
isomorphic in $\stmodcat{B}$ to a simple module $S'$ with $f\cdot
S'=0$;

$(2)$ For all simple $B$-modules $V$ with $f\cdot V=0$,
$\Phi^{-1}(V)$ is isomorphic in $\stmodcat{A}$ to a simple module
$V'$ with $e\cdot V'=0$.

\smallskip
\noindent If the stable equivalence $\Phi_1: \stmodcat{eAe}\ra
\stmodcat{fBf}$, induced from $\Phi$, lifts to a derived
equivalence between $eAe$ and $fBf$, then $\Phi$ lifts to an iterated almost $\nu$-stable
derived equivalence between $A$ and $B$. \label{ThmLift}
\end{Theo}

%Surprisingly, Theorem
%\ref{ThmLift} has the following consequence which provides us with a large variety of algebras,
%called Frobenius-finite algebras (see Section \ref{sect4.3} for
%definition), such that every stable equivalence
%of Morita type between these algebras can be lifted to a derived equivalence. Thus the
%class of Frobenius-finite algebras shares many common algebraic and
%numerical invariants of derived and stable equivalences. Of course,
%representation-finite algebras are Frobenius-finite algebras, but
%the converse is not true.

The contents of the paper is outlined as follows. In Section \ref{secPre}, we fix notation and collect some basic facts needed in our later proofs. In Section \ref{secPropStM}, we begin with a review of aspects on stable equivalences of Morita type, and then discuss the relationship between stable equivalences of Morita type over algebras and their Frobenius-parts which play a prominent role in our question mentioned above. In Sections \ref{secLift} and \ref{sect4.3}, we prove the main results, Theorem \ref{ThmLift} and Theorem \ref{ThmRepFin}, respectively. In Section \ref{sect6}, we illustrate the procedure of lifting stable equivalences of Morita type to derived equivalences discussed in the paper by two examples from the modular representation theory of finite groups. This shows that our results can be applied to verify the Brou\'e's abelian defect group conjecture for some cases. We end this section by a few open questions suggested by the main results in the paper.

\section{Preliminaries}\label{secPre}

In this section, we shall recall basic definitions and facts
required in our proofs.

\subsection{General notation on derived categories}
Throughout this paper, unless specified otherwise, all algebras will
be finite-dimensional algebras over a fixed field $k$. All modules will be
finitely generated unitary left modules.

Let $\cal C$ be an additive category. For two morphisms
$f:X\rightarrow Y$ and $g:Y\rightarrow Z$ in $\cal C$, the
composite of $f$ with $g$ is written as $fg$, which is a morphism
from $X$ to $Z$. But for two functors $F:\mathcal{C}\rightarrow
\mathcal{D}$ and $G:\mathcal{D}\rightarrow\mathcal{E}$ of
categories, their composite is denoted by $GF$. For an object $X$
in $\mathcal{C}$, we denote by $\add(X)$ the full subcategory of
$\cal C$ consisting of all direct summands of finite direct sums of
copies of $X$.

We denote by $\C{C}$ the category of complexes $\cpx{X}=(X^i,
d_X^i)$ over ${\cal C}$, where $X^i$ is an object in $\cal C$ and
the differential $d_X^i: X^i\ra X^{i+1}$ is a morphism in $\cal C$
with $d_X^id_X^{i+1}=0$ for each $i\in \mathbb{Z}$. The homotopy
category of $\mathcal C$ is denoted by $\K{\mathcal C}$. When $\cal C$ is
an abelian category, the derived category of ${\cal C}$ is denoted by
$\D{\cal C}$. The full subcategories of $\K{\cal C}$ and $\D{\cal C}$
consisting of bounded complexes over $\mathcal{C}$ are denoted by
$\Kb{\mathcal C}$ and $\Db{\mathcal C}$, respectively.

Let $A$ be an algebra. The category of all $A$-modules is denoted by
$\modcat{A}$; the full subcategory of $A$-mod consisting of
projective (respectively, injective) modules is denoted by
$\pmodcat{A}$ (respectively, $\imodcat{A}$).  $D$ is the usual
duality $\Hom_k(-, k)$. The duality $\Hom_A(-, A)$
from $\pmodcat{A}$ to $\pmodcat{A\opp}$ is denoted by $(-)^*$, that
is, for each projective $A$-module $P$, the projective
$A\opp$-module $\Hom_A(P, A)$ is denoted by $P^*$. We denote by
$\nu_A$ the Nakayama functor $D\Hom_A(-,A): A\mbox{-proj}\ra
A\mbox{-inj}$, which is an equivalence with $\nu^{-1}_A=\Hom_A(DA,-)$.  The stable module category $\stmodcat{A}$ of $A$ has
the same objects as $\modcat{A}$, and the morphism set $\StHom_A(X, Y)$
of two $A$-modules $X$ and $Y$ in $\stmodcat{A}$ is the quotient of $\Hom_A(X, Y)$
modulo the homomorphisms that factorize through projective modules.
As usual, we simply write $\Kb{A}$ and $\Db{A}$ for
$\Kb{\modcat{A}}$ and $\Db{\modcat{A}}$, respectively. It is well
known that $\Kb{A}$ and $\Db{A}$ are triangulated categories. For a
complex $\cpx{X}$ in $\K{A}$ or $\D{A}$, the complex $\cpx{X}[n]$ is
obtained from $\cpx{X}$ by shifting $\cpx{X}$ to the left by $n$
degrees.

\smallskip
For $X\in\modcat{A}$, we use $P(X)$ (respectively,
$I(X)$) to denote the projective cover (respectively, injective
envelope) of $X$. As usual, the syzygy and co-syzygy of $X$ are
denoted by $\Omega(X)$ and $\Omega^{-1}(X)$, respectively. The socle
and top, denoted by $\soc(X)$ and $\top(X)$, are the largest
semisimple submodule and the largest semisimple quotient module of
$X$, respectively.

\smallskip
A homomorphism $f: X\ra Y$ of $A$-modules is called a \emph{radical
map} if, for any module $Z$ and homomorphisms $h: Z\ra X$ and $g:
Y\ra Z$, the composite $hfg$ is not an isomorphism. A complex
over $\modcat{A}$ is called a \emph{radical} complex if all its
differential maps are radical maps. Every complex over $\modcat{A}$
is isomorphic in the homotopy category $\K{A}$ to a radical complex.
It is easy to see that if two radical complex $\cpx{X}$ and
$\cpx{Y}$  are isomorphic in $\K{A}$, then $\cpx{X}$ and $\cpx{Y}$
are isomorphic in $\C{A}$.

Two algebras $A$ and $B$ are said to be \emph{stably equivalent} if
their stable module categories $\stmodcat{A}$ and $\stmodcat{B}$ are
equivalent as $k$-categories, and \emph{derived equivalent} if
their derived categories $\Db{A}$ and $\Db{B}$ are equivalent as
triangulated categories. A triangle equivalence $F:\Db{A}\ra\Db{B}$
is called a {\em derived equivalence} between $A$ and $B$.

For derived equivalences, Rickard gave a nice characterization in
\cite{RickMoritaTh}. He showed that two algebras are derived
equivalent if and only if there is a complex $\cpx{T}$ in
$\Kb{\pmodcat{A}}$ satisfying

\smallskip
(1) $\Hom_{\Db{A}}(\cpx{T},\cpx{T}[n])=0$ for all  $n\ne 0$, and

(2) $\add(\cpx{T})$ generates $\Kb{\pmodcat{A}}$ as a triangulated
category

\smallskip
{\parindent=0pt such that $B\simeq\End_{\Kb{A}}(\cpx{T})$.}

A complex in $\Kb{\pmodcat{A}}$ satisfying the above two conditions
is called a {\em tilting complex} over $A$.
% By the above condition (2),
%each indecomposable projective $A$-module is a direct summand of
%$T^i$ for some integer $i$.
It is known that, given a derived equivalence $F$ between $A$ and
$B$, there is a unique (up to isomorphism) tilting complex $\cpx{T}$
over $A$ such that $F(\cpx{T})\simeq B$. This complex $\cpx{T}$ is
called a tilting complex \emph{associated} to $F$.

Recall that a complex $\cpx{\Delta}$ in $\Db{B\otimes_kA\opp}$ is
called a {\em two-sided tilting complex} provided that there is
another complex $\cpx{\Theta}$ in $\Db{A\otimes_kB\opp}$ such that
$\cpx{\Delta}\otimesL_{A}\cpx{\Theta}\simeq B$ in
$\Db{B\otimes_kB\opp}$ and
$\cpx{\Theta}\otimesL_{B}\cpx{\Delta}\simeq A$ in
$\Db{A\otimes_kA\opp}$. In this case, the functor
$\cpx{\Delta}\otimesL_{A}-: \Db{A}\ra \Db{B}$ is a derived
equivalence. A derived equivalence of this form is said to be {\em
standard}. For basic facts on the derived functor $-\otimesL -$, we
refer the reader to \cite{Weibel}.

\subsection{Almost $\nu$-stable derived equivalences\label{subsect2.2}}
In \cite{HuXi3}, a special kind of derived equivalences was
introduced, namely the almost $\nu$-stable derived equivalences.
Recall that a derived equivalence $F:\Db{A}\ra\Db{B}$ is called an
{\em almost $\nu$-stable derived equivalence} if the following two
conditions are satisfied:

\medskip
$(1)$ The tilting complex $\cpx{T}=(T^i, d^i)_{i\in \mathbb{Z}}$
associated to $F$ has only nonzero terms in negative degrees, that
is, $T^i=0$ for all $i>0$. In this case, the tilting complex
$\cpx{\bar{T}}$ associated to the quasi-inverse $G$ of $F$ has only
nonzero terms in positive degrees, that is, $\bar{T}^i=0$ for all
$i<0$ (see \cite[Lemma 2.1]{HuXi3}).

\medskip $(2)$
${\add(\bigoplus_{i<0}T^{i})=\add(\bigoplus_{i<0}\nu_AT^{i})}$ and
${\add(\bigoplus_{i>0}\bar{T}^i)=\add(\bigoplus_{i>0}\nu_B\bar{T}^i)}$.

\medskip
As was shown in \cite{HuXi3}, each almost $\nu$-stable derived
equivalence between $\Db{A}$ and $\Db{B}$ induces a stable equivalence
between $A$ and $B$. Thus $A$ and $B$ share many common invariants
of both derived and stable equivalences.

For the convenience of the reader, we briefly recall the construction in \cite{HuXi3}.

Suppose that $A$ and $B$ are two algebras and that
$F:\Db{A}\ra\Db{B}$ is a derived equivalence such that the tilting
complex associated to $F$ has no nonzero terms in positive degrees.
By \cite[Lemma 3.1]{HuXi3}, for each $X\in\modcat{A}$, one can fix a
radical complex $\cpx{\bar{Q}_X}\simeq F(X)$ in $\Db{B}$:
$$0\lra \bar{Q}_X^0\lra \bar{Q}_X^1\lra\cdots\lra \bar{Q}_X^n\lra 0$$
with $\bar{Q}_X^i$ projective for all $i>0$. Moreover, the complex
of this form is unique up to isomorphism in $\Cb{B}$. For $X, Y$ in
$\modcat{A}$, this induces an isomorphism
$$\phi: \Hom_A(X, Y)\lra \Hom_{\Db{B}}(\cpx{\bar{Q}_X}, \cpx{\bar{Q}_Y}).$$
Then a functor $\bar{F}: \stmodcat{A}\ra\stmodcat{B}$
, called the {\em stable functor} of $F$, was defined in \cite{HuXi3} as
follows: For each $X$ in $\modcat{A}$, we set
$$\bar{F}(X):=\bar{Q}_X^0.$$ For any morphism $f:
X\ra Y$ in $\modcat{A}$, we denote by $\underline{f}$ its image
in $\StHom_A(X, Y)$. By \cite[Lemma 2.2]{HuXi3}, the map $\phi(f)$
in $\Hom_{\Db{B}}(\cpx{\bar{Q}_X}, \cpx{\bar{Q}_Y})$ can be
presented by a chain map $\cpx{g}=(g^i)_{i\in \mathbb{Z}}$. Then we define $$\bar{F}:
\StHom_A(X, Y)\lra \StHom_B(\bar{F}(X), \bar{F}(Y)),\quad
\underline{f}\mapsto \underline{g}^0.$$ It was shown in \cite{HuXi3}
that $\bar{F}: \stmodcat{A}\ra\stmodcat{B}$ is indeed a
well-defined functor fitting into the following commutative diagram (up to
isomorphism)
$$\xy
  (0,0)*+{\stmodcat{B}}="c",
  (40,0)*+{\Db{B}/\Kb{\pmodcat{B}}}="d",
  (0,15)*+{\stmodcat{A}}="a",
  (40,15)*+{\Db{A}/\Kb{\pmodcat{A}}}="b",
  (80, 0)*+{\Db{B}}="dbb",
  (80, 15)*+{\Db{A}}="dba",
  {\ar^{\bar{F}} "a";"c"},
  {\ar^{F'} "b";"d"},
  {\ar^(.4){\Sigma_A}, "a";"b"},
  {\ar^(.4){\Sigma_B}, "c";"d"},
  {\ar_(.4){\rm can.} "dba";"b"},
  {\ar_(.4){\rm can.} "dbb";"d"},
  {\ar^{F} "dba";"dbb"},
\endxy$$
where $\Db{A}/\Kb{\pmodcat{A}}$ is a Verdier quotient, the functor
$\Sigma_A: \stmodcat{A}\ra \Db{A}/\Kb{\pmodcat{A}}$ is induced by
the canonical embedding $\modcat{A}\ra \Db{A}$, and $F'$ is the
triangle equivalence which is uniquely  determined (up to
isomorphism) by the commutative square in the right-hand side of the
above diagram.

One can easily check that, up to isomorphism, the stable functor
$\bar{F}$ is independent of the choices of the complexes
$\cpx{\bar{Q}_X}$. Moreover, if two derived equivalences are naturally
isomorphic, then so are their stable functors.

For a self-injective algebra $A$, it was shown in
$\cite{RickDstable}$ that the functor $\Sigma_{A}$ is a triangle
equivalence. Denote the composite
$$\Db{A}\lraf{\rm can.} \Db{A}/\Kb{\pmodcat{A}}\lraf{\Sigma_A^{-1}}\stmodcat{A}$$
by $\eta_{A}: \Db{A}\ra \stmodcat{A}$. Thus, if $A$ and $B$ are
self-injective algebras, then there is a uniquely determined (up to
isomorphism) equivalence functor $\Phi_{F}: \stmodcat{A}\ra
\stmodcat{B}$ such that the diagram
$$\xy
  (0,0)*+{\stmodcat{A}}="c",
  (25,0)*+{\stmodcat{B}}="d",
  (0,15)*+{\Db{A}}="a",
  (25,15)*+{\Db{B}}="b",
  {\ar^{\eta_A} "a";"c"},
  {\ar^{\eta_B} "b";"d"},
  {\ar^{F}, "a";"b"},
  {\ar^{{\Phi_{F}}}, "c";"d"},
\endxy$$
is commutative up to isomorphism. In this case, we say that the
stable equivalence $\Phi_F$ is {\em induced by the derived
equivalence $F$} or $\Phi_F$ {\em lifts to a derived equivalence}.

In general, a derived equivalence does not give rise to a stable
equivalence, nor the converse thereof. However, if a derived equivalence $F$ is almost
$\nu$-stable, then its stable functor $\bar{F}$ is a stable
equivalence \cite[Theorem 3.7]{HuXi3}. So we introduce the following definition:

If a stable equivalence $\Phi$ between arbitrary algebras is isomorphic to the stable functor $\bar{F}$ of an almost $\nu$-stable derived equivalence $F$, then we say that the
stable equivalence $\Phi$ is {\em induced by the almost
$\nu$-stable derived equivalence $F$}, or $\Phi$ {\em lifts to the almost $\nu$-stable derived  equivalence} $F$. If a stable equivalence
$\Phi$ can be written as a composite
$\Phi\simeq\Phi_1\circ\Phi_2\circ\cdots\circ\Phi_m$ of stable
equivalences with $\Phi_i$ or $\Phi_i^{-1}$ induced by an almost
$\nu$-stable derived equivalence for all $i$, then we say that
$\Phi$ is {\em induced by an iterated almost $\nu$-stable derived
equivalence} or $\Phi$ {\em lifts to an iterated almost $\nu$-stable derived equivalence} (see \cite{hu-iterated}).

Actually, the above two kinds of stable equivalences $\bar{F}$ and
$\Phi_F$ induced by derived equivalences are
compatible with each other when our consideration restricts to
self-injective algebras. Let $F:\Db{A}\ra\Db{B}$ is a derived
equivalence between two self-injective algebras. By the
above diagrams, if the tilting complex
associated to $F$ has no nonzero terms in positive degrees, then $F$
is an almost $\nu$-stable derived equivalence and the stable functor
$\bar{F}$ is isomorphic to the functor $\Phi_F$ defined above. If
the tilting complex $\cpx{T}$ associated to $F$ has nonzero terms in
positive degrees, then $F$ can be written as a composite $F\simeq
F_1\circ F_2^{-1}$ such that both $F_1$ and $F_2$ are almost
$\nu$-stable derived equivalences, and thus $\Phi_F\simeq
\Phi_{F_1}\circ \Phi_{F_2}^{-1}\simeq \bar{F}_1\circ\bar{F}_2^{-1}$.
Here we can take $F_2$ to be $[m]$ for which $\cpx{T}[-m]$ has no
nonzero terms in positive degrees. This shows that $\Phi_F$ lifts to an iterated almost $\nu$-stable derived equivalence.

Let us remark that if a derived equivalence $F$ is
standard and almost $\nu$-stable, then $\bar{F}$ is
a stable equivalence of Morita type (\cite[Theorem 5.3]{HuXi3}).
This is compatible with (and generalizes) the result \cite[Corollary
5.5]{RickDFun} of Rickard which says that $\Phi_F$ is a stable
equivalence of Morita type provided that $F$ is a standard derived
equivalence between two self-injective algebras.

\subsection{Frobenius parts and $\nu$-stable idempotent elements\label{sect2.3}}

In this subsection, we recall the definition of the Frobenius part of an algebra, which was introduced in \cite{mv} and related to the
Nakayama functor, and collect some basic facts related to idempotent elements.

Let $A$ be an algebra, and let $e$ be an idempotent element in $A$. It is well known that
$Ae\otimes_{eAe}-:\modcat{eAe}\ra\modcat{A}$ is a full embedding
and induces a full embedding
$$\lambda: \stmodcat{eAe}\lra \stmodcat{A}.$$
There is another functor $eA\otimes_A-: \modcat{A}\ra \modcat{eAe}$, such that the functors $Ae\otimes_{eAe}-$ and $eA\otimes_A-$ induce mutually
inverse equivalences between $\add(Ae)$ and $\pmodcat{eAe}$. Further, the functor $eA\otimes_A-$ induces a triangle equivalence between the homotopy categories $\Kb{\add(Ae)}$ and $\Kb{\pmodcat{eAe}}$. In particular, if $P\in \add(Ae)$, then $Ae\otimes_{eAe}eA\otimes_AP\simeq P$ as $A$-modules.
Moreover, we have the following facts.

\begin{Lem} Let $A$ be an algebra, and let $e$ be an idempotent element in $A$. For
a simple $A$-module $S$, we define $\Delta_e(S):=Ae\otimes_{eAe}eS$
and denote by $P(S)$ the projective cover of $S$.
Suppose that $eS\neq 0$. Then

 $(1)$ $\Delta_e(S)$ is isomorphic to a quotient module of $P(S)$ and $e\cdot\rad(\Delta_e(S))=0$;

 $(2)$ If $e\cdot\rad(P(S))\neq 0$, then $\Delta_e(S)$ is non-projective.
 \label{LemIdemProp0}
\end{Lem}

{\it Proof.} $(1)$ Applying $Ae\otimes_{eAe}eA\otimes_A-$ to the
epimorphism $P(S)\ra S$, we get an epimorphism $Ae\otimes_{eAe}
eP(S)\ra \Delta_e(S)$. Since $eS\neq 0$, the projective cover $P(S)$
of $S$ is in $\add(Ae)$, and therefore $Ae\otimes_{eAe}e P(S)\simeq
P(S)$ by the equivalence between $\add(Ae)$ and $\pmodcat{eAe}$.
Hence $\Delta_e(S)$ is isomorphic to a quotient module of $P(S)$.
Thus $\Delta_e(S)$ has $S$ as a single top. Applying $eA\otimes_A-$
to the short exact sequence $0\ra \rad(\Delta_e(S))\ra
\Delta_e(S)\ra S\ra 0$, we have a short exact sequence
$$0\lra e\cdot\rad(\Delta_e(S))\lra e\cdot\Delta_e(S)\lraf{h} eS\lra 0.$$
The middle term $e\cdot\Delta_e(S)\simeq eAe\otimes_{eAe}eS\simeq
eS$. This implies that $h$ must be an isomorphism and
$e\cdot\rad(\Delta_e(S))=0$.

$(2)$ Suppose contrarily that $\Delta_e(S)$ is projective. Then the
epimorphism $P(S)\ra \Delta_e(S)$ splits. This forces that
$\Delta_e(S)\simeq P(S)$. By assumption, we have
$e\cdot\rad(P(S))\neq 0$, while $e\cdot\rad(\Delta_e(S))=0$. This is
a contradiction. $\square$

\medskip
We say that an idempotent element $e$ in $A$ is  {\em $\nu$-stable}
provided that $\add(\nu_AAe)=\add(Ae)$. That is, for each
indecomposable direct summand $P$ of $Ae$, the corresponding
injective module $\nu_AP$ is still a direct summand of $Ae$. Clearly, the module $Ae$ is projective-injective.
Note that the notion of $\nu$-stable idempotents is left-right symmetric, although it is defined by using left modules.
In fact, $\add(\nu_AAe)=\add(Ae)$ if and only if $\add(eA)=\add\big(\nu_{A\opp}(eA)\big)$ because $D(\nu_AAe)\simeq DD(eA)\simeq eA$ and $D(Ae)\simeq \nu_{A\opp}(eA)$. Moreover, we have the following lemma.

\begin{Lem}
Let $A$ be an algebra, and let $e$ be a $\nu$-stable idempotent in $A$. Then

$(1)$  $\add(\top(Ae))=\add(\soc(Ae))$.

$(2)$ If $\add(Ae)\cap\add\big(A(1-e)\big)=\{0\}$, then $\soc(eA)$ is an ideal of $A$. Moreover, $\soc(Ae)=\soc(eA)$.
\label{lem2.2a}
\end{Lem}

{\it Proof.} (1) Since $\top(Ae)=\soc(\nu_AAe)$, the statement (1) follows from the definition of $\nu$-stable idempotents.

(2) By our assumption, it follows from \cite[Section 9.2]{DrozdKir} that $\soc(Ae)$ is an ideal of $A$. It follows from (1) that
$(1-e)\soc(Ae)=0$. Thus $\soc(Ae)=\big((1-e)\cdot\soc(Ae)\big)\oplus \big(e\cdot\soc(Ae)\big)=e\cdot\soc(Ae)\subseteq eA$. Moreover, for each $r\in\rad(A)$, the left $A$-module homomorphism $\phi_r: A\ra A, x\mapsto xr$ is a radical map. The restriction of $\phi_r$ to any indecomposable direct summand $X$ of $Ae$ cannot be injective. Otherwise, $\phi_r|_X$ is split since $X$ is injective, and $\phi_r$ is not a radical map. This is a contradiction. Hence $\soc(X)\subseteq\Ker\,\phi_r$, and $\soc(Ae)\subseteq \Ker\,\phi_r$. This means that $\soc(Ae)\cdot r=0$. Consequently $\soc(Ae)\subseteq\soc(eA)$. The duality $\Hom_A(-, A)$ takes $Ae$ to $eA$, and $A(1-e)$ to $(1-e)A$. This implies that $\add(eA)\cap\add\big((1-e)A\big)=\{0\}$. Similarly, we have $\soc(eA)\subseteq\soc(Ae)$, and therefore $\soc(eA)=\soc(Ae)$. $\square$

\medskip
A projective $A$-module $P$ is called {\em $\nu$-stably projective} if $\nu_A^iP$ is projective for all $i>0$.
We denote by $\stp{A}$ the full subcategory of $\pmodcat{A}$ consisting of all $\nu$-stably projective $A$-modules.
Clearly, $\stp{A}$ is closed under taking direct summands and finite direct sums. The two notions of $\nu$-stable idempotents and $\nu$-stably projective modules are closely related. Actually we have the following lemma.

\begin{Lem}
  Let $A$ be an algebra. Then the following hold.

  $(1)$ If $e$ is a $\nu$-stable idempotent in $A$, then $\add(Ae)\subseteq\stp{A}$.

  $(2)$ If $e$ is an idempotent in $A$ such that $\add(Ae)=\stp{A}$, then $e$ is $\nu$-stable.

  $(3)$ There is an $\nu$-stable idempotent $e$ in $A$ such that $\add(Ae)=\stp{A}$.

  $(4)$ All the modules in $\stp{A}$ are projective-injective.
  \label{lem2.2}
\end{Lem}

{\it proof.} (1) Let $P\in \add(Ae)$. Then, by definition, the module $\nu_AP\in\add(\nu_AAe)=\add(Ae)$, and consequently $\nu_A^iP$ belongs to $\add(Ae)$ for all $i>0$. Hence $P$ is a $\nu$-stably projective $A$-module, that is, $P\in\stp{A}$.

(2) Since $Ae\in\stp{A}$, the $A$-module $\nu_A^i(Ae)$ is projective for all $i>0$. This further implies that $\nu_AAe$ is projective and  $\nu_A^i(\nu_AAe)$ is projective for all $i>0$. Hence $\nu_AAe\in\stp{A}=\add(Ae)$, and $\add(\nu_AAe)\subseteq\add(Ae)$. Since $\nu_A$ is an equivalence from $\pmodcat{A}$ to $\imodcat{A}$, the categories $\add(\nu_AAe)$ and $\add(Ae)$ have the same number of isomorphism classes of indecomposable objects. Hence $\add(\nu_AAe)=\add(Ae)$, that is, the idempotent $e$ is $\nu$-stable.

$(3)$  Since $\stp{A}$ is a full subcategory of $\pmodcat{A}$, there is an idempotent $e$ in $A$ such that $\add(Ae)=\stp{A}$. The statement (3) then follows from (2).

$(4)$ By definition, all the modules in $\stp{A}$ are projective. By $(3)$, there is a $\nu$-stable idempotent $e$ such that $\add(Ae)=\stp{A}$. This implies that all the modules in $\stp{A}$ are in $\add(Ae)=\add(\nu_AAe)$, and  subsequently they are also injective. $\square$

\medskip

If $e$
is an idempotent element in $A$ such that $\add(Ae)=\stp{A}$, then the algebra $eAe$ is called the \emph{Frobenius part} of
$A$, or the {\em associated self-injective algebra} of $A$. Clearly, the Frobenius part of $A$ is unique up to Morita equivalence.

\begin{Lem} Let $A$ be an algebra, and let $e$ be an idempotent element of $A$.
Then we have the following:

$(1)$ For $Y\in\add(Ae)$ and $X\in\modcat{A}$, there is an
isomorphism induced by the functor $$eA\otimes_A-: \Hom_A(Y, X)\lra \Hom_{eAe}(eY, eX).$$

$(2)$ There is a natural isomorphism $e(\nu_AY)\simeq \nu_{eAe}(eY)$ for all $Y\in\add(Ae)$.
%  $$\Hom_A(X, \nu_AP)\lra \Hom_{eAe}(e(X), e(\nu_{A}P)).$$

$(3)$ If $e$ is $\nu$-stable, then $eAe$ is a self-injective algebra.

$(4)$ Suppose that $e$ is $\nu$-stable. If the algebra $A$ has no
semisimple direct summands, then neither does the algebra $eAe$.
\label{LemIdemProp}
\end{Lem}

{\it Proof.}
(1) is well known (see, for example, \cite[Proposition
2.1, p33]{AusReiten}).

(2) follows from $(1)$ and the following isomorphisms
$$\begin{array}{rl}
\nu_{eAe}(eY) & = D\Hom_{eAe}(eY, eAe)\\
 & \simeq D\Hom_A(Y, Ae)\simeq D(Y^*\otimes_AAe)\\
 &\simeq\Hom_A(Ae,
 D(Y^*)) \simeq e(\nu_AY)
\end{array}$$

(3) follows immediately from (2) (see also \cite{mv}).

(4) Since the functor $eA\otimes_A-: \add(Ae)\ra \pmodcat{eAe}$ is
an equivalence, each indecomposable projective $eAe$-module is isomorphic to $eY$
for some indecomposable $A$-module $Y$ in $\add(Ae)$. By definition,
we have $\add(Ae)=\add(\nu_AAe)$, which means that $Y$ is
projective-injective and $\soc(Y)\in\add(\top(Ae))$. Since $A$ has
no semisimple direct summands, the module $Y$ is not simple. Thus
$Y$ has at least two composition factors in $\add(\top(Ae))$ and
consequently $eY$ has at least two composition factors. Hence $eY$
is not simple. This implies that the algebra $eAe$ has no semisimple direct
summands. $\square$

\medskip
The following lemma is easy. But, for the convenience of the reader, we include here a proof.

\begin{Lem}
Let $A$ be an algebra, and let $M$ be an $A$-module which is a
generator for $\modcat{A}$, that is, $\add(_AA)\subseteq \add(M)$.
Suppose that $X$ is an $A$-module. Then $\Hom_A(M, X)$ is a
projective $\End_A(M)$-module if and only if $X\in\add(M)$.
  \label{lemEndMproj}
\end{Lem}

{\it Proof.} Clearly, if $X\in \add(M)$, then $\Hom_A(M,X)$ is a projective $\End_A(M)$-module. Now, suppose that $\Hom_A(M,X)$ is projective for an $A$-module $X$.  Without loss of generality, we may assume that $A$ is a
basic algebra. Then ${}_AA$ is a direct summand of $M$, that is,
$M\simeq A\oplus N$ for some $A$-module $N$. Since $\Hom_A(M, X)$ is
a projective $\End_A(M)$-module, there is some $M_X\in\add(M)$ such
that $\Hom_A(M, M_X)\simeq \Hom_A(M, X)$ as $\End_A(M)$-modules. By
Yoneda isomorphism, there is an $A$-module homomorphism $f:M_X\ra X$
such that $\Hom_A(M, f)$ is an isomorphism, that is, $\Hom_A(A,
f)\oplus \Hom_A(N, f)$ is an isomorphism. This implies that
$\Hom_A({}_AA, f)$ is an isomorphism, and therefore so is $f$. $\square$

\medskip
Finally, we point out the following elementary facts on
Nakayama functors.

(1) For any $A$-module $M$ and projective $A$-module $P'$, there is a
natural isomorphism: $D\Hom_A(P',M)\simeq \Hom_A(M,\nu_AP')$. More general, for any $\cpx{P}\in \Kb{\pmodcat{A}}$ and $\cpx{X}\in \Kb{A}$, there is is an isomorphism of $k$-spaces: $D\Hom_{\Kb{A}}(\cpx{P},\cpx{X})\simeq \Hom_{\Kb{A}}(\cpx{X},\nu_A\cpx{P})$.

(2) Let $M$ be a fixed generator for $A$-mod, and let
$\Lambda:=\End_A(M)$. Then, for each projective $A$-module $P'$,
there is a natural isomorphism
$\nu_{\Lambda}\Hom_A(M, P')\simeq \Hom_A(M, \nu_AP').$

\section{Stable equivalences of Morita
type}\label{secPropStM}

As a special kind of stable equivalences, Brou\'e introduced the notion of stable equivalences of Morita type (see, for example, \cite{Broue1994}), which is a combination of Morita and stable equivalences. In this section, we shall first collect some basic properties of
stable equivalences of Morita type, and then give conditions for
lifting stable equivalences of Morita type to Morita equivalences
which are, of course, special kinds of derived equivalences. The results in this
section will be used in Section \ref{secLift} for the proof of the
main result, Theorem \ref{ThmLift}.

\subsection{Basic facts on stable equivalences of Morita type\label{sect3.1}}
Let $A$ and $B$ be two $k$-algebras over a field $k$.
Following \cite{Broue1994}, we say that two bimodules $_AM_B$ and
${}_BN_A$ define a \emph{stable equivalence of Morita type} between $A$ and
$B$ if the following conditions hold:

\medskip
(1) The one-sided modules $_AM, M_{B}, {}_BN$ and $N_{A}$ are all
projective;

(2) $M\otimes_BN\simeq A\oplus P$ as $A$-$A$-bimodules for some
projective $A$-$A$-bimodule $P$, and $N\otimes_AM\simeq B\oplus Q$
as $B$-$B$-bimodules for some projective $B$-$B$-bimodule $Q$.

\medskip
In this case, we have two exact functors
$T_M=M\otimes_B-: \modcat{B} \ra \modcat{A}$ and  $T_N ={}_BN\otimes_A-:
\modcat{A}\ra\modcat{B}$. Analogously, the bimodules $P$ and $Q$ define
two exact functors $T_P$ and $T_Q$, respectively. Note that the
images of $T_P$ and $T_Q$ consist of projective modules. Moreover,
the functor $T_N$ induces an equivalence $\Phi_N:
\stmodcat{A}\ra\stmodcat{B}$ for stable moduule categories. The functor $\Phi_N$ is called a {\em
stable equivalence of Morita type}. Similarly, we have $\Phi_M$
which is a quasi-inverse of $\Phi_N$.

Clearly, $P=0$ if and only if $Q=0$. In this situation,
we come back to the notion of Morita equivalences.

It would be interesting to replace the word ``projective" by ``flat" or ``Gorenstein flat or projective" in the above definition and to deduce the corresponding ``stable" theory. We refrain  from these considerations here.

\medskip
For stable equivalences of Morita type, we have the following basic
facts.

\begin{Lem}
Let $A$ and $B$ be algebras without
semisimple direct summands. Suppose that $_AM_B$ and $_BN_A$ are two
bimodules without projective direct summands and define a stable
equivalence of Morita type between $A$ and $B$. Write $_AM\otimes_BN_A\simeq A\oplus
  P$ and $_BN\otimes_AM_B\simeq B\oplus Q$ as bimodules. Then the
  following hold.

\smallskip
$(1)$ $(M\otimes_B-, N\otimes_A-)$ and $(N\otimes_A-,
     M\otimes_B-)$ are adjoint pairs of functors.

$(2)$ $\add(\nu_AP)=\add({}_AP)$ and $\add(\nu_BQ)=\add({}_BQ)$.

$(3)$ $N\otimes_AP\in\add(_BQ)$, and $M\otimes_BQ\in\add(_AP)$.

$(4)$ For each indecomposable $A$-module $X\not\in\add(_AP)$, the
$B$-module $N\otimes_AX$ is the direct sum of an indecomposable module
$\bar{X}\not\in\add({}_BQ)$ and a module $X'\in\add(_BQ)$.

$(5)$ If $S$ is a simple $A$-module with $\Hom_A({}_AP, S)=0$,
 then $N\otimes_AS$ is simple with $\Hom_B({}_BQ, N\otimes_AS)=0$.

$(6)$ Suppose that $A/\rad(A)$ and $B/\rad(B)$ are separable.
 If  $S$ is a simple $A$-module with $\Hom_A({}_AP, S)\neq 0$, then
 $N\otimes_AS$ is indecomposable and non-simple with  both
 $\soc(N\otimes_AS)$ and $\top(N\otimes_AS)$ in $\add(\top({}_BQ))$.
\label{LemPropStM}
\end{Lem}

{\it Proof}. (1) This follows from \cite{DugasVilla} and
\cite{LiuSummands} (see also \cite[Lemma 4.1]{ChenPanXi}).

(2) For an $A$-module $X$, we see that $P\otimes_AX$ is in
$\add(_AP)$. In fact, if we take a surjective homomorphism
$(_AA)^n\ra X$, then we get a surjective map $P\otimes_AA^n\ra
P\otimes_AX$. Since $_AP\otimes_AX$ is projective for all
$A$-modules $X$, we know that $P\otimes_AX$ is a direct summand of
$_AP^n$.

We have the following isomorphisms
$$\begin{array}{rl}
  \nu_B(N\otimes_AX) & =D\Hom_B(N\otimes_AX, B)\\
                     & \simeq D\Hom_A(X, M\otimes_BB)   \quad ( \mbox{by } (1))\\
                     &\simeq D\Hom_A(X, A\otimes_AM)\\
                     & \simeq D(\Hom_A(X, A)\otimes_AM)  \quad (\mbox{because }{}_AM \mbox{ is projective})\\
                     & \simeq \Hom_A(M, \nu_AX)   \quad ( \mbox{ by adjointness })\\
                     & \simeq \Hom_B(B, N\otimes_A\nu_AX)   \quad ( \mbox{by (1)})\\
                     & \simeq N\otimes_A(\nu_AX).
\end{array}$$
Similarly, for a $B$-module $Y$, we have
$\nu_A(M\otimes_BY)\simeq M\otimes_B(\nu_BY)$. Thus
$\nu_A(M\otimes_BN\otimes_AA)\simeq
M\otimes_BN\otimes_A(\nu_AA)$, and consequently $\nu_AA\oplus\nu_AP\simeq (A\oplus
P)\otimes_A(\nu_AA)$. Hence $\nu_AP\simeq
P\otimes_A(\nu_AA)\in\add(_AP)$, and therefore
$\add(_AP)\subseteq \add(\nu_AP)$. Since $\nu_A$ is an equivalence from
$\pmodcat{A}$ to $A$-{inj},
we deduce that $\add(_AP)= \add(\nu_AP)$ just by counting the number of
indecomposable direct summands of $_AP$ and $\nu_AP$. Similarly, we have
$\add(_BQ)=\add(\nu_BQ)$. This proves $(2)$.

(3) It follows from $N\otimes_A(A\oplus
P)\simeq N\otimes_AM\otimes_BN\simeq (B\oplus
Q)\otimes_BN$ that $N\otimes_AP\simeq
Q\otimes_BN$ as bimodules. In particular, as a left $B$-module, $N\otimes_AP$ is isomorphic to $Q\otimes_BN$ which is in
$\add(_BQ)$.  Hence $N\otimes_AP\in\add(_BQ)$. Similarly, $M\otimes_BQ\in\add(_AP)$.

(4) Suppose that $X$ is an indecomposable $A$-module and $X\not\in
\add(_AP)$. Let $N\otimes_AX=\bar{X}\oplus X'$ be a decomposition of
$N\otimes_AX$ such that $\bar{X}$ has no direct summands in
$\add(_BQ)$ and $X'\in\add(_BQ)$. If $\bar{X}=0$, then
$N\otimes_AX\in\add(_BQ)$ and consequently $X\oplus
P\otimes_AX\simeq M\otimes_B(N\otimes_AX)\in\add(_AP)$ by (3). This
is a contradiction. Hence $\bar{X}\neq 0$. Suppose that $\bar{X}$
decomposes, say $\bar{X}=Y_1\oplus Y_2$ with $Y_i\neq 0$ for $i=1,
2$. Clearly, $M\otimes_BY_i\not\in \add(_AP)$ for $i=1,2.$ It
follows that both $M\otimes_BY_1$ and $M\otimes_BY_2 $ have
indecomposable direct summands which are not in $\add({}_AP)$.
However, we have an isomorphism $X\oplus P\otimes_AX\simeq
M\otimes_BN\otimes_AX\simeq M\otimes_BY_1\oplus M\otimes_BY_2\oplus
M\otimes_BX'$, and $X$ is the only indecomposable direct summand of
$X\oplus P\otimes_AX$ not in $\add({}_AP)$.  But $X$ is the only
indecomposable direct summand of $X\oplus P\otimes_AX$ with
$X\not\in \add(_AP)$. This contradiction shows that $\bar{X}$ must
be indecomposable.

(5) By (1) and \cite[Lemma 3.2]{Xi2} together with the proof of \cite[Lemma 4.5]{Xi2}
we have $P\simeq P^*$ as
$A$-$A$-bimodules. Note that this was proved in \cite[Proposition 3.4]{DugasVilla} with additional conditions that $A/\rad(A)$ and $B/\rad(B)$ are separable. If $\Hom_A(P, S)=0$, then $P\otimes_AS\simeq
P^*\otimes_AS\simeq\Hom_A(P, S)=0$. Thus, we have
$M\otimes_BN\otimes_AS\simeq S\oplus P\otimes_AS=S$. Note that
$N\otimes_A-$ is an exact and faithful functor since
$A_A\in\add(N_A)$. We denote by $\ell(X)$ the length of the
composition series of $X$. It follows that $\ell(N\otimes_AX)\geq
\ell(X)$ for all $A$-modules $X$. Similarly, $\ell(M\otimes_BY)\geq
\ell(Y)$ for all $B$-modules $Y$. Consequently, we have
$1=\ell(S)\leq \ell(N\otimes_AS)\leq
\ell(M\otimes_BN\otimes_AS)=\ell(S)=1$. This implies that
$N\otimes_AS$ is a simple $B$-module. Finally, $\Hom_B({}_BQ,
N\otimes_AS)\simeq \Hom_A(M\otimes_BQ, S)=0$ by (1) and (3).

(6) Let $e$ be an idempotent element in $A$ such that
$\add(_AAe)=\add({}_AP)$ and $\add(Ae)\cap\add\big(A(1-e)\big)=\{0\}$, and let $f$ be
an idempotent element in $B$ such that $\add({}_BBf)=\add({}_BQ)$ and $\add(Bf)\cap\add\big(B(1-f)\big)=\{0\}$.
Then $e$ and $f$ are $\nu$-stable idempotents, and the modules $eA_A$ and ${}_BBf$ are projective-injective. Consequently, the $B$-$A$-bimodule $Bf\otimes_keA$ is also projective-injective and $\add\big((B\otimes_kA)(f\otimes e)\big)\cap\add\big((B\otimes_kA)(1-f\otimes e)\big)=\{0\}$. By Lemma \ref{lem2.2a}, $\soc(eA_A)$, $\soc({}_BBf)$ and $\soc(Bf\otimes_keA)$ are ideals of $A$, $B$ and $B\otimes_kA\opp$, respectively, and $\soc(Ae)=\soc(eA)$.  Since the algebras
$A/\rad(A)$ and
$B/\rad(B)$ are separable, we have $\soc({}_BBf\otimes_k
eA_A)=\soc(Bf)\otimes_k\soc(eA)$. By assumption, the bimodule $N$ has no projective
direct summands. Particularly, $N$ has no direct summands in
$\add(Bf\otimes_keA)$. This is equivalent to that
$\soc(Bf\otimes_keA)N=0$ by \cite[Section 9.2]{DrozdKir}. That is, $\soc(Bf)N\soc(eA)=0$. As $N_A$ is projective, we have $N\otimes_A\soc(eA)\simeq N\soc(eA)$. Thus
$$\soc(Bf)(N\otimes_A\soc(eA))\simeq \soc(Bf)(N\soc(eA))= \soc(Bf)N\soc(eA)=0.$$
This means that the $B$-module $N\otimes_A\soc(eA)$ has no direct
summands in $\add({}_BQ)$. Now let $S$ be a simple $A$-module with
$\Hom_A(P, S)\neq 0$. Then $S$ is in $\add(\top(Ae))=\add(\soc(Ae))$.
Since $\soc(Ae)=\soc(eA)$, we have $S\in\add({}_A\soc(eA))$, and consequently the $B$-module $N\otimes_AS$
has no direct summands in $\add({}_BQ)$. Now, by (4), the
module $N\otimes_AS$ is indecomposable. Suppose that $N\otimes_AS$
is simple. Then $M\otimes_B(N\otimes_AS)$ must be indecomposable by
the above discussion. However, we have an isomorphism
$M\otimes_B(N\otimes_AS)\simeq S\oplus P\otimes_AS$. This forces
that $P\otimes_AS=0$ and implies that $\Hom_A(P,S)\simeq \Hom_A(_AP,_AA)\otimes_AS\simeq P\otimes_AS=0$, a contradiction. Hence
$N\otimes_AS$ is an indecomposable non-simple $B$-module. Since
$\Hom_A({}_AP, S)\neq 0$, there is a sequence
$P\lraf{f}S\lraf{g}\nu_AP$ with $f$ surjective and $g$ injective.
Applying the exact functor $N\otimes_A-$, we get a new sequence
$N\otimes_AP\lraf{N\otimes_Af}N\otimes_AS\lraf{N\otimes_Ag}N\otimes_A\nu_AP$
with $N\otimes_Af$ surjective and $N\otimes_Ag$ injective. By (2)
and (3), we see that both $\soc(N\otimes_AS)$ and $\top(N\otimes_AS)$ are
in $\add(\top({}_BQ))$.
$\hfill\square$

\medskip
Now, let us make a few comments on the
separability condition in the above
lemma. Suppose that $A$ is a finite-dimensional $k$-algebra over a field $k$. That $A/\rad(A)$ is a separable algebra over $k$ is
equivalent to that the center of $\End_A(S)$ is a separable
extension of $k$ for each simple $A$-module $S$. Thus, if
$A$ satisfies the separability condition, then so do its quotient
algebras and the algebras of the form $eAe$ with $e$ an idempotent
element in $A$. The separability condition seems not to be a strong restriction and can be satisfied actually by many interesting classes of algebras. Here, we mention a few: A finite-dimensional $k$-algebra $A$
satisfies the separability condition if one of the following is fulfilled:

$\bullet$ $k$ is a perfect field. For example, a finite field, an algebraically closed field, or a field of characteristic zero.

$\bullet$ $A$ is given by a quiver with relations.

$\bullet$ $A$ is the group algebra $kG$ of a finite group $G$ (see,
for example, \cite[Lemma 1.28, p. 183]{RepFinGroup}).

\subsection{Stable equivalences of Morita type at different levels\label{sect3.2}}
We say that a stable equivalence $\Phi: \stmodcat{A}\ra
\stmodcat{B}$ of Morita type {\em lifts to a Morita equivalence} if
there is a Morita equivalence $F: \modcat{A}\ra \modcat{B}$ such
that the diagram
$$\xy
  (0,0)*+{\stmodcat{A}}="c",
  (25,0)*+{\stmodcat{B}}="d",
  (0,15)*+{\modcat{A}}="a",
  (25,15)*+{\modcat{B}}="b",
  {\ar^{\rm can.} "a";"c"},
  {\ar^{\rm can.} "b";"d"},
  {\ar^{F}, "a";"b"},
  {\ar^{{\Phi}}, "c";"d"},
\endxy$$
of functors is commutative up to isomorphism, where the vertical functors are the canonical ones.

The following proposition tells us when a stable equivalence of
Morita type lifts to a Morita equivalence.

\begin{Prop} Let $A$ and $B$ be algebras without
semisimple direct summands. Suppose that $_AM_B$ and $_BN_A$ are two
bimodules without projective direct summands and define a stable
equivalence of Morita type between $A$ and $B$. Write $_AM\otimes_BN_A\simeq A\oplus
P$ and $_BN\otimes_AM_B\simeq B\oplus Q$ as bimodules. Then the following are equivalent:

\smallskip
 $(1)$ $N\otimes_A-: \modcat{A}\ra \modcat{B}$ is an equivalence, that is, $P=0=Q$.

 $(2)$ $N\otimes_AS$ is a simple $B$-module for every simple $A$-module $S$.

\smallskip
{\parindent=0pt If} $A/\rad(A)$ and $B/\rad(B)$ are separable, then
the above statements are further equivalent to the following two
equivalent conditions:

$(3)$ The stable equivalence $\Phi_N$ induced by $N\otimes_A-$ lifts
to a Morita equivalence.

$(4)$ $N\otimes_AS$ is isomorphic in $\stmodcat{B}$ to a simple
$B$-module for every simple $A$-module $S$.
 \label{PropMorEqCond}
\end{Prop}

{\it Proof.}
$(1)\Rightarrow(2)$ is trivial, since $N\otimes_A-$ is a Morita equivalence in case that $P=0=Q$.

$(2)\Rightarrow(1)$ was proved by Liu \cite{LiuSimpletoSimple} under
the condition that the ground field $k$ is splitting for
both $A$ and $B$. Here, we give a proof which is independent of the
ground field. Suppose that $P\neq 0$. Let $\{S_1, \cdots, S_m\}$ be
a complete set of non-isomorphic simple $A$-modules in
$\add(\top({}_AP))$. Then, since $_AP$ is a projective-injective module and $A$ has no semisimple direct summands, the
indecomposable direct summands of $_{A}P$  cannot be simple, and consequently all $S_i$ are not projective and $S_i\not\in \add({}_AP)$.  Thus it follows from $S_i\oplus P\otimes_AS_i\simeq M\otimes_BN\otimes_AS_i$ that
the simple $B$-modules $N\otimes_AS_i$ and $N\otimes_AS_j$ are not
isomorphic whenever $i\neq j$. Using the adjoint pair in Lemma
\ref{LemPropStM} (1), we get the following isomorphisms:
$$\begin{array}{rl} \End_A(S_i)\oplus \Hom_A(P\otimes_AS_i, \bigoplus_{j=1}^m S_j)
&\simeq \Hom_A(S_i\oplus P\otimes_AS_i, \bigoplus_{j=1}^m S_j) \\
         &\simeq \Hom_A(M\otimes_BN\otimes_AS_i, \bigoplus_{j=1}^m S_j) \\
         &\simeq \Hom_B(N\otimes_AS_i, \bigoplus_{j=1}^m N\otimes_AS_j) \\
        &\simeq \End_B(N\otimes_AS_i) \\
         &\simeq \underline{\End}_B(N\otimes_AS_i)\quad (N\otimes_AS_i \mbox{ is a non-projective simple module }) \\ & \simeq  \underline{\End}_A(S_i) \\ & \simeq \End_A(S_i),
      \end{array}$$
which implies that $\Hom_A(P\otimes_AS_i, \bigoplus_{j=1}^m S_j)=0$.
However, the $A$-module $P\otimes_AS_i$ belongs to $\add(_AP)$ and is nonzero since
$P\otimes_AS_i\simeq P^*\otimes_A S_i\simeq \Hom_A({}_AP, S_i)\neq
0$. This implies that $\Hom_A(P\otimes_AS_i, \bigoplus_{j=1}^m
S_j)\neq 0$, a contradiction. Thus $P=0$, and therefore $Q=0$.

Note that $(1)\Rightarrow (3)\Rightarrow (4)$ is obvious.

For the rest of the proof, we assume that $A/\rad(A)$ and $B/\rad(B)$ are separable.

It remains to prove ``$(4)\Rightarrow (2)$". According to Lemma
\ref{LemPropStM} (5), this can be done by showing
$\Hom_A({}_AP, S)=0$ for all simple $A$-modules $S$. Let $S$ be an
arbitrary simple $A$-module. If $S$ is not projective, then it follows from
Lemma \ref{LemPropStM} (6) that $\Hom_A(P, S)=0$ since $\Phi_N(S)$
isomorphic to a simple $B$-module in $\stmodcat{B}$. If $S$ is
projective, then it cannot be in $\add({}_AP)$. Otherwise, $S$ is
projective-injective and $A$ has a semisimple block, contradicting
to our assumption. Hence $\Hom_A({}_AP, S)=0$. $\square$

\medskip
Now, we recall the restriction procedure of stable equivalences of
Morita type from \cite[Theorem 1.2]{ChenPanXi}. Suppose that $A$ and
$B$ are two algebras without semisimple direct summands, and that
$_AM_B, {}_BN_A$ are two bimodules without projective direct
summands, and define a stable equivalence of Morita type between $A$ and
$B$. If $e$ and $f$ are idempotent elements in $A$ and $B$, respectively, such
that $M\otimes_BNe\in\add(Ae)$ and $\add(Bf)=\add(Ne)$, then the
bimodules $eMf$ and $fNe$ define a stable equivalence of Morita type
between $eAe$ and $fBf$, that is, the diagram
$$\xy
  (0,0)*+{\stmodcat{eAe}}="c",
  (35,0)*+{\stmodcat{fBf}}="d",
  (0,15)*+{\stmodcat{A}}="a",
  (35,15)*+{\stmodcat{B}}="b",
  {\ar^{\lambda} "c";"a"},
  {\ar^{\lambda} "d";"b"},
  {\ar^{\Phi_N}, "a";"b"},
  {\ar^{\Phi_{fNe}}, "c";"d"},
\endxy$$
is commutative up to isomorphism, where $\lambda$ is defined in Section \ref{sect2.3}.

The following lemma describes the restriction of
stable equivalences in terms of simple modules.

\begin{Lem} Let $A$ and $B$ be algebras without
semisimple direct summands such that $A/\rad(A)$ and $B/\rad(B)$ are
separable. Suppose that $e$ and $f$ are idempotent
elements in $A$ and $B$, respectively. Let  $\Phi: \stmodcat{A}\ra\stmodcat{B}$ be a stable
equivalence of Morita type such that

\smallskip
$(1)$ For each simple $A$-module $S$ with $e\cdot S=0$, the
$B$-module $\Phi(S)$ is isomorphic in $\stmodcat{B}$ to a simple
module $S'$ with $f\cdot S'=0$;

$(2)$ For each simple $B$-module $T$ with $f\cdot T=0$, the
$A$-module $\Phi^{-1}(T)$ is isomorphic in $\stmodcat{A}$ to a
simple module $T'$ with $e\cdot T'=0$.

\smallskip
{\parindent=0pt Then} there is, up to isomorphism, a unique stable
equivalence $\Phi_1:\stmodcat{eAe}\ra\stmodcat{fBf}$ of Morita type
such that the following diagram of functors
$$\xy
  (0,0)*+{\stmodcat{eAe}}="c",
  (25,0)*+{\stmodcat{fBf}}="d",
  (0,15)*+{\stmodcat{A}}="a",
  (25,15)*+{\stmodcat{B}}="b",
  {\ar^{\lambda} "c";"a"},
  {\ar^{\lambda} "d";"b"},
  {\ar^{\Phi}, "a";"b"},
  {\ar^{\Phi_1}, "c";"d"},
\endxy$$
is commutative (up to isomorphism).  \label{LemRes}
\end{Lem}

{\it Proof.}
We may assume that the stable equivalence $\Phi$ of Morita
type between $A$ and $B$ is defined by bimodules $_AM_B$ and $_BN_A$ without
projective direct summands, that is, $\Phi\simeq \Phi_N$ which is induced by
the functor $_BN\otimes_A-$. Suppose that $M\otimes_BN\simeq A\oplus
P$ and $N\otimes_AM\simeq B\oplus Q$ as bimodules. By the assumption (1)
and Lemma \ref{LemPropStM} (6), we have $\Hom_A({}_AP, S)=0$ for all
simple $A$-modules $S$ with $e\cdot S=0$. This implies that
${}_AP\in\add(Ae)$, and consequently $M\otimes_BNe\simeq Ae\oplus
Pe\in\add(Ae)$. Now, for each simple $B$-module $T$ with $f\cdot
T=0$, it follows from the assumption (2) that $\Hom_A(Ae,
M\otimes_BT)=0$. This is equivalent to $\Hom_B(N\otimes_AAe, T)=0$
 by Lemma \ref{LemPropStM} (1). Hence $Ne\simeq N\otimes_AAe\in
\add(Bf)$. Similarly, we get $_BQ\in\add(Bf)$ and
$M\otimes_BBf\in\add(Ae)$, and consequently $Bf$ is in
$\add(N\otimes_AM\otimes Bf)\subseteq\add(N\otimes_AAe)=\add(Ne)$.
Therefore $\add(Ne)=\add(Bf)$. Using \cite[Theorem 1.2]{ChenPanXi},
we get the desired commutative diagram. The functor $\Phi_1$ is
uniquely determined up to natural isomorphism because $\lambda$ is a
full embedding.
$\square$

\medskip
The next proposition shows that a stable equivalence of Morita type between
$A$ and $B$ may lift to a Morita equivalence provided that certain
`restricted' stable equivalence lifts to a Morita equivalence.

\medskip

\begin{Prop}
  Let $A$ and $B$ be two algebras without semisimple direct
  summands such that $A/\rad(A)$ and $B/\rad(B)$ are separable, and let $e$ and $f$ be
  idempotents in $A$ and $B$, respectively.  Suppose that there is a commutative (up to isomorphism) diagram
$$\xy
  (0,0)*+{\stmodcat{eAe}}="c",
  (25,0)*+{\stmodcat{fBf}}="d",
  (0,15)*+{\stmodcat{A}}="a",
  (25,15)*+{\stmodcat{B}}="b",
  {\ar^{\lambda} "c";"a"},
  {\ar^{\lambda} "d";"b"},
  {\ar^{\Phi}, "a";"b"},
  {\ar^{\Phi_1}, "c";"d"},
\endxy$$
with $\Phi$ and $\Phi_1$ being stable equivalences of Morita type,
and that the following conditions hold:

\smallskip
  $(1)$  For each simple $A$-module $S$ with $e\cdot S=0$, the $B$-module $\Phi(S)$ is isomorphic in $\stmodcat{B}$ to a simple $B$-module.

  $(2)$ For each  simple $B$-module $T$ with $f\cdot T=0$, the $A$-module $\Phi^{-1}(T)$ is isomorphic in $\stmodcat{A}$ to a simple $A$-module.

\noindent If $\Phi_1$ lifts to a Morita equivalence, then $\Phi$ lifts to a Morita
equivalence. \label{PropLiftMorita}
\end{Prop}

{\it Proof.}
Suppose that $_AM_B$ and $_BN_A$ are two bimodules without projective direct summands and define a
stable equivalence of Morita type between $A$ and $B$ such that
$\Phi$ is induced by $N\otimes_A-$. Assume that $M\otimes_BN\simeq A\oplus
P$ and $N\otimes_AM\simeq B\oplus Q$ as bimodules. We shall prove $P=0$.

Assume contrarily that $P\neq 0$. Let $S$ be a simple $A$-module
with $\Hom_A({}_AP, S)\neq 0$. Then $S$ cannot be projective,
Otherwise, $S$ would be a direct summand of ${}_AP$ which is
projective-injective, and $A$ would have a semisimple direct
summand. We shall prove that $N\otimes_AS$ is isomorphic to a simple
$B$-module $T$, and this will lead to a contradiction by Lemma
\ref{LemPropStM} (6).

First, we claim that $eS\neq 0$. Otherwise, it would follow from the assumption (1) that $\Phi(S)$ is isomorphic to a simple $B$-module, a contradiction by Lemma \ref{LemPropStM} (6). Hence $eS\neq 0$, $P(S)\in\add(Ae)$. This implies also that each indecomposable direct summand of $P$ is in $\add(Ae)$ since we can choose a simple module $S$ for each of such summands so that $\Hom_A(P,S)\ne 0$. Consequently, we have ${}_AP\in\add(Ae)$. Similarly we have ${}_BQ\in\add(Bf)$. Since $\Phi_1$ lifts to a Morita equivalence, the module $\Phi_1(eS)$ is isomorphic in
$\stmodcat{fBf}$ to a simple $fBf$-module $fT$ with $T$ a simple
$B$-module. Set $\Delta_e(S):=Ae\otimes_{eAe}eS$ and $\Delta_f(T):=Bf\otimes_{fBf}fT$. By the given
commutative diagram, we get an isomorphism in $\stmodcat{B}$
$$(*) \quad  N\otimes_A\Delta_e(S)\simeq \Delta_f(T).$$
Now, we claim that $N\otimes_A\Delta_e(S)$ and $\Delta_f(T)$ are actually isomorphic in $\modcat{B}$. To prove this, it suffices to show that
$N\otimes_A\Delta_e(S)$ is indecomposable and non-projective.

Note that $P(S)$ is an indecomposable non-simple projective-injective module. In fact, it follows from $\Hom_A(P,S)\ne 0$ that $P(S)$ is a direct summand of $P$ which is projective-injective by Lemma \ref{LemPropStM}(2). Moreover, $P(S)$ is non-simple because $A$ has no semisimple direct summands. Thus, we have $\soc(P(S))\subseteq\rad(P(S))$.
Since
$\add(\nu_AP)=\add({}_AP)$ by Lemma \ref{LemPropStM} (2), we have
$\add(\top({}_AP))=\add(\soc({}_AP))$. Hence $\soc(P(S))\in\add(\top({}_AP))\subseteq\add(\top(Ae))$. Consequently $e\cdot\soc(P(S))\neq 0$, and $e\cdot\rad(P(S))\neq 0$. By Lemma \ref{LemIdemProp0}, the $A$-module $\Delta_e(S)$, which is a quotient module of $P(S)$, is not projective. This implies that $N\otimes_A\Delta_e(S)$ is not projective.

By Lemma \ref{LemPropStM} (4), to prove that $N\otimes_A\Delta_e(S)$ is indecomposable, we have to show that $N\otimes_A\Delta_e(S)$ has no direct summands in $\add({}_BQ)$. Suppose contrarily that this is false and $Q_1\in\add({}_BQ)$ is an indecomposable direct summand of $N\otimes_A\Delta_e(S)$. We consider the exact sequence
$$(**)\quad 0\lra N\otimes_A\rad (\Delta_e(S))\lra N\otimes_A\Delta_e(S)\lra N\otimes_AS\lra 0.$$
Then $\Hom_A(N\otimes_A\rad (\Delta_e(S)), Q_1)\neq 0$. Otherwise it follows from the exact sequence $(**)$ that $Q_1$ has
to be a direct summand of $N\otimes_A S$ which is indecomposable by Lemma \ref{LemPropStM} (6). Thus $Q_1\simeq N\otimes_AS$. However, since
$S$ is not projective, the module $N\otimes_AS$ cannot be projective. This leads to a contradiction. Thanks to the formula $\Hom_A(\nu^{-1}_AY,X)\simeq D\Hom_A(X,Y)$ for any $A$-module $X$ and any projective $A$-module $Y$, we have
$$\Hom_A\big(\nu_A^{-1}(M\otimes_BQ_1), \rad(\Delta_e(S))\big)\simeq D\Hom_A(\rad(\Delta_e(S)),
M\otimes_BQ_1)\simeq D\Hom_B(N\otimes_A\rad(\Delta_e(S)), Q_1)\neq
0.$$ By Lemma \ref{LemPropStM} (2) and (3), we know that $\nu^{-1}_A(M\otimes_BQ_1)\in \add(P)$, $\Hom_A({}_AP,
\rad(\Delta_e(S)))\neq 0$, and $e\cdot \rad(\Delta_e(S))\simeq
\Hom_A(Ae, \rad(\Delta_e(S)))\neq 0$. This contradicts to Lemma \ref{LemIdemProp0} and shows that the $B$-module
$N\otimes_A\Delta(S)$ is indecomposable.

Hence $N\otimes_A\Delta_e(S)\simeq \Delta_f(T)$ in $\modcat{B}$. Together with the exact sequence $(**)$ above, we deduce that $N\otimes_AS$ is isomorphic to a quotient module of $\Delta_f(T)$. By Lemma \ref{LemPropStM} (6), the socle of $N\otimes_AS$ is in $\add(\top({}_BQ))$. Since ${}_BQ\in\add(Bf)$, we have $\soc(N\otimes_AS)\in\add(\top(Bf))$. However, it follows from Lemma \ref{LemIdemProp0} that $\Delta_f(T)$ has top $T$ and $f\cdot\rad(\Delta_f(T))=0$. This means that $\rad(\Delta_f(T))$ has no composition factors in $\add(\top(Bf))$, and $T$ is the only quotient module of $\Delta_f(T)$ with the socle in $\add(\top(Bf))$. This yields that $N\otimes_AS\simeq T$, which contradicts to Lemma \ref{LemPropStM} (6). Hence $P=0$. This implies that $N\otimes_A-$ is a Morita
equivalence between the module categories of the algebras $A$ and
$B$. $\square$

\section{From stable equivalences of Morita type to derived
equivalences}\label{secLift}

In this section, we shall prove the main result, Theorem \ref{ThmLift}.
We first make some preparations.

\subsection{Extending derived
equivalences}\label{secExtTilt}

Let $A$ be an algebra over a field $k$, and let $e$ be a
$\nu$-stable idempotent element in $A$. In this subsection, we shall
show that a tilting complex over $eAe$ can be extended to an
tilting complex over $A$ which defines an almost $\nu$-stable derived equivalence.

First, we fix some terminology on approximations.

Let $\cal C$ be a category, $\cal D$ be a full subcategory
of $\cal C$, and $X$ be an object in $\cal C$. A morphism $f: D\ra X$
in $\cal C$ is called a \emph{right} $\cal D$-\emph{approximation}
of $X$ if $D\in {\cal D}$ and the induced map Hom$_{\cal C}(-,f)$:
Hom$_{\cal C}(D',D)\ra$ Hom$_{\cal C}(D',X)$ is surjective for
every object $D'\in {\cal D}$.  A morphism $f:X\ra Y$ in $\cal C$
is called \emph{right minimal} if any morphism $g: X\ra X$ with
$gf=f$ is an automorphism. A minimal right $\cal D$-approximation of
$X$ is a right $\cal D$-approximation of $X$, which is right
minimal. Dually, there is the notion of a \emph{left} $\cal
D$-\emph{approximation} and a \emph{minimal left} $\cal
D$-\emph{approximation}. The subcategory $\mathcal D$ is said to be \emph{functorially finite} in $\mathcal C$ if every object in $\mathcal C$ has a right and left $\mathcal D$-approximation.

\medskip
The following proposition shows that if an idempotent element $e$ in $A$ is $\nu$-stable,
then every tilting complex over $eAe$
can be extended to a tilting complex over $A$. This result extends \cite[Theorem 4.11]{miy} in which
algebras are assumed to be symmetric.

\medskip
\begin{Prop}
 Let $A$ be an arbitrary algebra over a field $k$, and let
 $e$ be a $\nu$-stable idempotent element in $A$. Suppose that $\cpx{Q}$ is
 a complex in $\Kb{\add(Ae)}$ with $Q^i=0$ for all $i>0$ such that $e\cpx{Q}$ is a tilting
 complex over $eAe$. Then there is a complex $\cpx{P}$ of $A$-modules such that
 $\cpx{Q}\oplus\cpx{P}$ is a tilting complex over $A$ and induces an
 almost $\nu$-stable derived equivalence between $A$ and the endomorphism algebra of the tilting module. \label{PropETilt}
\end{Prop}

{\it Proof}. For convenience, we shall abbreviate $\Hom_{\Kb{A}}(-, -)$ to $\Hom(-, -)$ in the proof.
Assume that $\cpx{Q}$ is of the following form:
$$0\lra Q^{-n}\lra \cdots\lra Q^{-1}\lra Q^0\lra 0$$
for some fixed natural number $n$. Note that $\add(\cpx{Q})$
is a functorially finite subcategory in $\Kb{A}$ since both $\Hom(\cpx{Q},
\cpx{X})$ and $\Hom(\cpx{X}, \cpx{Q})$ are finite-dimensional for
each $\cpx{X}\in\Kb{A}$. Thus, there is a minimal right
$\add(\cpx{Q})$-approximation $f_n: \cpx{Q_n}\ra A[n]$. The
following construction is standard. Let $\cpx{P_n}:=A[n]$. We define inductively a complex $\cpx{P_i}$ for each $i\le n$ by taking the following  distinguished triangle in $\Kb{\pmodcat{A}}$
$$(\star)\qquad \cpx{P_{i-1}}\lra \cpx{Q_i}\lraf{f_i}\cpx{P_i}\lra \cpx{P_{i-1}}[1],$$
where $f_i$ is a minimal right $\add(\cpx{Q})$-approximation of $\cpx{P_i}$ and where
$\cpx{P_{i-1}}[1]$ is a radical complex isomorphic in
$\Kb{\pmodcat{A}}$ to the  mapping cone of $f_i$. In the following,
we shall prove that $\cpx{Q}\oplus\cpx{P_0}$ is a tilting complex over $A$ and
induces an almost $\nu$-stable derived equivalence.

Clearly, by definition, $\add(\cpx{Q}\oplus \cpx{P_0})$ generates $\Kb{\pmodcat{A}}$. It remains to show that
$\Hom(\cpx{Q}\oplus \cpx{P_0},\cpx{Q}[m]\oplus \cpx{P_0}[m])=0$
for all $m\ne 0$. We shall prove this in four steps.

\medskip
(a) We show that $\Hom(\cpx{Q},\cpx{Q}[m])=0$ for all $m\ne 0$.

\medskip
In fact, it follows from the equivalence $eA\otimes_A-: \add(Ae)\ra
\add(\pmodcat{eAe})$ that $eA\otimes_A-$ induces a triangle
equivalence $\Kb{\add(Ae)}\ra \Kb{\pmodcat{eAe}}$. Since $e\cpx{Q}$
is a tilting complex over $eAe$, we see that
$\Hom(e\cpx{Q},e\cpx{Q}[m])=0$ for all $m\ne 0$. Therefore, for the
complex $\cpx{Q}\in \Kb{\add(Ae)}$, we have
$\Hom(\cpx{Q},\cpx{Q}[m])=0$ for all $m\ne 0$.

\medskip
(b) We claim that $\Hom(\cpx{Q},\cpx{P_0}[m])=0$ for all $m\ne 0$.

\medskip
Indeed, applying
$\Hom(\cpx{Q}, -)$ to the above triangle $(\star)$, we obtain a long exact sequence
$$\cdots\lra \Hom(\cpx{Q}, \cpx{P_{i-1}}[m])\lra \Hom(\cpx{Q},\cpx{Q_i}[m])\lra
\Hom(\cpx{Q}, \cpx{P_{i}}[m])\lra \Hom(\cpx{Q},
\cpx{P_{i-1}}[m+1])\lra\cdots $$
for each integer
$i\leq n$. Since $\Hom(\cpx{Q},\cpx{Q}[m])=0$ for all $m\ne 0$, one gets
$$\Hom(\cpx{Q}, \cpx{P_{i-1}}[m])\simeq \Hom(\cpx{Q}, \cpx{P_{i}}[m-1])$$
for all $m<0$. Thus, for all $m<0$, we have
 $$\Hom(\cpx{Q}, \cpx{P_0}[m])\simeq \Hom(\cpx{Q}, \cpx{P_1}[m-1])\simeq\cdots\simeq\Hom(\cpx{Q},\cpx{P_n}[m-n])\simeq\Hom(\cpx{Q},A[m])=0.$$

To prove that $\Hom(\cpx{Q}, \cpx{P_0}[m])=0$ for $m>0$, we shall show by induction on $i$ that $\Hom(\cpx{Q},\cpx{P_i}[m])=0$ for all $m>0$ and all $i\le n$.

Clearly, for $i=n$, we have $\Hom(\cpx{Q},\cpx{P_n}[m])=0$ for all $m>0$.
Now,  assume inductively that $\Hom(\cpx{Q},\cpx{P_j}[m])=0$ for all $m>0$ and all $i\le j\le n$.
Since $f_i$ is a right $\add(\cpx{Q})$-approximation of $\cpx{P_i}$, the induced map
$\Hom(\cpx{Q}, f_i)$ is surjective. Thus $\Hom(\cpx{Q},\cpx{P_{i-1}}[1])=0$ by (a). The long exact sequence, together with (a) and the induction hypothesis, yields that
$\Hom(\cpx{Q},\cpx{P_{i-1}}[m])=0$ for all $m>1$. Thus $\Hom(\cpx{Q},\cpx{P_i}[m])=0$ for all $m>0$ and all $i\le n$. In particular, $\Hom(\cpx{Q},\cpx{P_0}[m])=0$ for all $m>0$. This completes the proof of (b).

\medskip
(c) $\Hom(\cpx{P_0},\cpx{Q}[m])=0$ for all $m\ne 0$.

\medskip
To prove (c),  let $\Delta$ be the endomorphism algebra of $e\cpx{Q}$. Since
$\add(Ae)=\add(\nu_AAe)$, the algebra $eAe$ is a self-injective
algebra by Lemma \ref{LemIdemProp} (3). Thanks to \cite[Theorem 2.1]{Al-Rickard}, we see that $\Delta$ is also
self-injective. Let $G:\Db{eAe}\ra\Db{\Delta}$ be the derived
equivalence induced by the tilting complex $e\cpx{Q}$. Then $G(e\cpx{Q})$ is
isomorphic to $\Delta$. Since $\Delta$ is self-injective, we have
$\add(\nu_{\Delta}\Delta)=\add(_{\Delta}\Delta)$, and consequently
$\add(e\cpx{Q})=\add(\nu_{eAe}e\cpx{Q})$, or equivalently $\add(\cpx{Q})=\add(\nu_A\cpx{Q})$. Therefore
$\Hom(\cpx{P_0}, \cpx{Q}[m])\simeq D\Hom(\nu_A^{-1}\cpx{Q},
\cpx{P_0}[-m])=0$ for all $m\neq 0$.

\medskip
(d) Finally, we show that $\Hom(\cpx{P_0},\cpx{P_0}[m])=0$ for all $m\ne 0$.

\medskip
Indeed, we know that $G(eAe)$ is isomorphic to a complex $\cpx{V}$ in
$\Kb{\pmodcat{\Delta}}$ with $V^i=0$ for all $i<0$ (see, for instance,
\cite[Lemma 2.1]{HuXi3}) and have shown in (b) that
 $\Hom(\cpx{Q}, \cpx{P_0}[m])$ = $0$
for all $m\neq 0$. It follows that $\Hom_{\Kb{eAe}}\big(e\cpx{Q}$,
$e\cpx{P_0}[m]\big) = 0$ for all $m\neq 0$, and consequently
$G(e(\cpx{P_0}))$ is isomorphic in $\Db{\Delta}$ to a
$\Delta$-module. Thus
$$\begin{array}{rl}
\Hom(Ae, \cpx{P_0}[m]) & \simeq\Hom_{\Kb{eAe}}(eAe,
e(\cpx{P_0})[m])\\
  & \simeq \Hom_{\Db{eAe}}(eAe,e(\cpx{P_0})[m])\\
&
\simeq\Hom_{\Db{\Delta}}(\cpx{V},G(e(\cpx{P_0}))[m])\\ & =0\end{array}$$
for all $m>0$. By the construction of $\cpx{P_0}$, all the terms of
$\cpx{P_0}$ in non-zero degrees lie in $\add(Ae)$. Since $\cpx{P_0}$
is a radical complex, the term $P_0^m$ is zero  for all $m>0$. Otherwise we would have $\Hom(Ae, \cpx{P_0}[t])\neq 0$
for the maximal positive integer $t$ with $P_0^t\neq 0$.

Applying the functor $\Hom(\cpx{P_0}, -)$ to the triangle $(\star)$, we have
an exact sequence (for all $m$ and $i\leq n$)
$$\Hom(\cpx{P_0}, \cpx{Q_i}[m-1])\lra \Hom(\cpx{P_0}, \cpx{P_i}[m-1])\lra \Hom(\cpx{P_0}, \cpx{P_{i-1}}[m])\lra \Hom(\cpx{P_0}, \cpx{Q_i}[m]).$$
If $m<0$, then $\Hom(\cpx{P_0},
\cpx{Q_i}[m-1])=0=\Hom(\cpx{P_0}, \cpx{Q_i}[m])$, and
$\Hom(\cpx{P_0}, \cpx{P_i}[m-1])\simeq \Hom(\cpx{P_0},
\cpx{P_{i-1}}[m])$. Thus, for $m<0$, we get
$$ \Hom(\cpx{P_0}, \cpx{P_0}[m])\simeq \Hom(\cpx{P_0}, \cpx{P_1}[m-1])\simeq\cdots\simeq\Hom(\cpx{P_0}, \cpx{P_n}[m-n])=\Hom(\cpx{P_0}, A[m])=0.$$
Now we apply $\Hom(-, \cpx{P_0})$ to the triangle $(\star)$ and get an exact sequence (for all  $m$ and
$i\leq n$)
$$\Hom(\cpx{Q_i}, \cpx{P_0}[m])\lra\Hom(\cpx{P_{i-1}}, \cpx{P_0}[m])\lra\Hom(\cpx{P_i}, \cpx{P_0}[m+1])\lra\Hom(\cpx{Q_i}, \cpx{P_0}[m+1]).$$
If $m>0$, then $\Hom(\cpx{Q_i}, \cpx{P_0}[m])=0=\Hom(\cpx{Q_i},
\cpx{P_0}[m+1])$, and consequently $\Hom(\cpx{P_{i-1}},
\cpx{P_0}[m])\simeq\Hom(\cpx{P_i}, \cpx{P_0}[m+1])$. Thus, for
$m>0$, we have
$$\Hom(\cpx{P_0}, \cpx{P_0}[m])\simeq\Hom(\cpx{P_1}, \cpx{P_0}[m+1])\simeq\cdots\simeq\Hom(\cpx{P_n}, \cpx{P_0}[m+n])=\Hom(A, \cpx{P_0}[m])=0.$$
So, we have proved that $\cpx{T}:=\cpx{Q}\oplus\cpx{P_0}$ is a
tilting complex over $A$. Let $B$ be the endomorphism algebra of
$\cpx{T}$ and let $F: \Db{A}\ra\Db{B}$ be the derived equivalence
induced by $\cpx{T}$. Then $F(\cpx{Q})$ is isomorphic in $\Db{B}$ to
the $B$-module $\Hom(\cpx{T}, \cpx{Q})$ with the property that
$\add(\nu_B(\Hom(\cpx{T}, \cpx{Q})))=\add(\Hom(\cpx{T},\cpx{Q}))$,
since $\add(\cpx{Q})=\add(\nu_A\cpx{Q})$ and $F$ commutes with the
Nakayama functor (see \cite[Lemma 2.3]{HuXi3}). By the definition of
$\cpx{P_0}$, we infer that $F(A)$ is isomorphic to a complex with
terms in $\add (\Hom(\cpx{T}, \cpx{Q}))$ for all positive degrees,
and zero for all negative degrees. Thus, by \cite[Proposition
3.8]{HuXi3}, the derived equivalence $F$ is almost $\nu$-stable. If
we define $\cpx{P}:=\cpx{P_0}$, then Proposition \ref{PropETilt}
follows. $\square$

\medskip
{\it Remark. } In Proposition \ref{PropETilt}, if we replace the condition
``$Q^i=0$ for all $i>0$" by the dual condition ``$Q^i=0$ for all $i<0$", then a dual
construction gives us a tilting complex $\cpx{Q}\oplus\cpx{P}$, which induces the quasi-inverse of
an almost $\nu$-stable derived equivalence.

%(2) Proposition \ref{PropETilt} extends \cite[Theorem 4.11]{miy} in which
%algebras are assumed to be symmetric.

\iffalse
 Let $A$ be a finite-dimensional algebra and let $e$
be an $\nu$-stable idempotent. If $\cpx{Q}$ is a complex in
$\Kb{\add(Ae)}$ with $Q^i=0$ for all $i<0$ such that $e(\cpx{Q})$ is
a tilting complex over $eAe$, then, by a dual construction, there is
a complex $\cpx{P}$ such that $\cpx{Q}\oplus\cpx{P}$ is a tilting
complex over $A$ and the quasi-inverse of the derived equivalence,
which is induced by $\cpx{Q}\oplus\cpx{P}$, is almost $\nu$-stable.
\fi

\medskip
\begin{Lem}
Keep the assumptions and notation in {\rm Proposition \ref{PropETilt}}. Let
$B:=\End_{\Kb{\pmodcat{A}}}(\cpx{Q}\oplus \cpx{P})$, and let $f$ be the idempotent element
in $B$ corresponding to the summand $\cpx{Q}$. Then there is a commutative (up
to isomorphism) diagram of functors
  $$\xy
      (0,0)*+{\stmodcat{eAe}}="eae",
      (0,15)*+{\stmodcat{A}}="a",
      (25,0)*+{\stmodcat{fBf}}="fbf",
      (25,15)*+{\stmodcat{B}}="b",
      {\ar^{\Phi} "a";"b"},
      {\ar^{\Phi_1} "eae";"fbf"},
      {\ar^{\lambda} "eae";"a"},
      {\ar^{\lambda} "fbf";"b"}
  \endxy$$
such that

$(1)$ $\Phi$ is a stable equivalence of Morita type induced by an
almost $\nu$-stable derived equivalence.

$(2)$ $\Phi_1$ is a stable equivalence of Morita type induced by a
derived equivalence $G$ with $G(e\cpx{Q})\simeq fBf$.

$(3)$ $\Phi(S)$ is isomorphic in $\stmodcat{B}$ to a simple
$B$-module $S'$ with $f\cdot S'=0$ for all simple $A$-modules $S$
with $e\cdot S=0$.

$(4)$ $\Phi^{-1}(T)$ is isomorphic in $\stmodcat{A}$ to a simple
$A$-module $T'$ with $e\cdot T'=0$ for all simple $B$-modules $T$
with $f\cdot T=0$.
 \label{LemTiltStM}
\end{Lem}

{\it Proof.} We first show the existence of the commutative diagram
of functors and the statements (1) and (2). By Proposition
\ref{PropETilt}, there is a derived equivalence $F:\Db{A}\ra \Db{B}$
such that $F(\cpx{Q}\oplus\cpx{P})\simeq B$ and $F(\cpx{Q})\simeq
Bf$. Moreover, $F$ is an almost $\nu$-stable derived equivalence.
Since $e\cpx{Q}$ is a tilting complex over $eAe$, we know that
$\add(e\cpx{Q})$ generates $\Kb{\pmodcat{eAe}}$ as a triangulated
category. Equivalently, $\add(\cpx{Q})$ generates $\Kb{\add(Ae)}$ as
a triangulated category. Thus, the functor $F$ induces a triangle
equivalence between $\Kb{\add(Ae)}$ and $\Kb{\add(Bf)}$.

By \cite[Corollary 3.5]{RickDFun}, there is a standard derived equivalence which agrees with $F$ on $\Kb{\pmodcat{A}}$. So, we can assume that $F$ itself is a standard derived
equivalence, that is, there are complexes $\cpx{\Delta}\in\Db{B\otimes_kA\opp}$ and $\cpx{\Theta}\in\Db{A\otimes_kB\opp}$ such that
$\cpx{\Delta}\otimesL_A\cpx{\Theta}\simeq {}_BB_B$, $\cpx{\Theta}\otimesL_B\cpx{\Delta}\simeq {}_AA_A$ and $F=\cpx{\Delta}\otimesL_A-$. By \cite[Lemma
5.2]{HuXi3}, the complex $\cpx{\Delta}$ can be assumed as the following form
$$0\lra \Delta^0\lra \Delta^1\lra \cdots\lra\Delta^n\lra 0$$
such that $\Delta^i\in\add(Bf\otimes_keA)$ for all $i>0$ and $\Delta^0$ is projective as left and right modules, and that $\cpx{\Theta}$ can be chosen to equal $\HomP_B(\cpx{\Delta}, {}_BB)$. Moreover, we have $\cpx{\Delta}\otimesP_A\cpx{\Theta}\simeq {}_BB_B$ in $\Kb{B\otimes_kB\opp}$ and $\cpx{\Theta}\otimesP_B\cpx{\Delta}\simeq {}_AA_A$ in $\Kb{A\otimes_kA\opp}$.

Since all the terms of $\cpx{\Delta}$ are projective as right $A$-modules, it follows
that $F(\cpx{X})\simeq \cpx{\Delta}\otimesL_A\cpx{X}\simeq\cpx{\Delta}\otimesP_A\cpx{X}$ for all $\cpx{X}\in\Db{A}$. Hence
$\cpx{\Delta}\otimesP_AAe\simeq F(Ae)$ is isomorphic in $\Db{B}$ to a complex in $\Kb{Bf}$.
Moreover, for each $i>0$, the term $\Delta^i\otimes_AAe$ is in $\add(Bf)$
since $\Delta^i\in\add(Bf\otimes_keA)$. Thus $\Delta^0\otimes_AAe\in\add(Bf)$.
Hence
$\Delta^i\otimes_AAe\in\add(Bf)$
  for all integers $i$, and consequently all
   the terms of the complex $f\cpx{\Delta}e$:
   $$0\lra f\Delta^0e\lra f\Delta^1e\lra\cdots\lra f\Delta^ne\lra 0$$
   are projective as left $fBf$-modules.

Next, we show that $f\Delta^ie$ is projective as right $eAe$-modules for all $i$. We have the following isomorphisms in $\K{A\opp}$:
$$\begin{array}{rl}
    f\cpx{\Delta} &\simeq fB\otimesP_B\cpx{\Delta}\\
              & \simeq \Hom_B(Bf, {}_BB)\otimesP_B\cpx{\Delta}\\
              & \simeq \HomP_B(F(\cpx{Q}), {}_BB)\otimesP_B\cpx{\Delta}\\
              & \simeq \HomP_B(\cpx{\Delta}\otimesP_A\cpx{Q}, {}_BB)\otimesP_B\cpx{\Delta}\\
              & \simeq \HomP_A(\cpx{Q}, \HomP_B(\cpx{\Delta},{}_BB))\otimesP_B\cpx{\Delta}\\
              & \simeq \HomP_A(\cpx{Q}, \HomP_B(\cpx{\Delta},{}_BB)\otimesP_B\cpx{\Delta})\\
              & \simeq \HomP_A(\cpx{Q}, \cpx{\Theta}\otimesP_B\cpx{\Delta})\\
              & \simeq \HomP_A(\cpx{Q}, {}_AA).\\
\end{array}$$
Since $\cpx{Q}\in\Kb{\add(Ae)}$, the complex $\HomP_A(\cpx{Q}, {}_AA)$ is in $\Kb{\add(eA)}$. For each $i>0$, it follows from the fact that $\Delta^i\in\add(Bf\otimes_keA)$ that $f\Delta^i\in\add(eA)$. Thus, using the above isomorphism in $\K{A\opp}$, we see that $f\Delta^0$ is again in $\add(eA)$, and consequently $f\Delta^i\in\add(eA)$ for all $i$. Hence $f\Delta^ie$ is projective as right $eAe$-modules for all $i$.

Now we have
the following isomorphisms in $\Db{fBf\otimes_kfBf\opp}$:
$$\begin{array}{rl}
     f\cpx{\Delta}e\otimesL_{eAe}e\cpx{\Theta}f &\simeq
     f\cpx{\Delta}e\otimesP_{eAe}e\cpx{\Theta}f\\
     & \simeq
     (fB\otimesP_B\cpx{\Delta}\otimesP_AAe)\otimesP_{eAe}(eA\otimesP_A\cpx{\Theta}\otimesP_BBf)\\
     & \simeq
     fB\otimesP_B\cpx{\Delta}\otimesP_A\big(Ae\otimesP_{eAe}eA\otimesP_A(\cpx{\Theta}\otimesP_BBf)\big)\\
     &\simeq
     fB\otimesP_B\cpx{\Delta}\otimesP_A\cpx{\Theta}\otimesP_BBf \quad (\mbox{ because } \cpx{\Theta}\otimesP_BBf\in\Kb{\add(Ae)})\\
     &\simeq
     fB\otimes_BB\otimes_BBf\\ & \simeq fBf.
   \end{array}$$
Similarly, $e\cpx{\Theta}f\otimesL_{fBf}f\cpx{\Delta}e\simeq eAe$ in
$\Db{eAe\otimes_keAe\opp}$. Thus $f\cpx{\Delta}e$ is a two-sided
tilting complex and $f\cpx{\Delta}e\otimesL_{eAe}-:
\Db{eAe}\ra\Db{fBf}$ is a derived equivalence. Note that we have the following isomorphisms in $\Db{fBf}$:
$$f\cpx{\Delta}e\otimesL_{eAe}e\cpx{Q}\simeq
f\cpx{\Delta}e\otimesP_{eAe}e\cpx{Q}\simeq
f\cpx{\Delta}\otimesP_A\cpx{Q}\simeq fBf.$$ This means that $e\cpx{Q}$ is an associated tilting
complex of the functor $G:=f\cpx{\Delta}e\otimesL_{eAe}-$.

Since $F=\cpx{\Delta}\otimesP_A-$ is an almost $\nu$-stable, standard derived equivalence, it follows from
\cite[Theorem 5.3]{HuXi3} that $\Delta^0\otimes_A-$ induces a stable equivalence of Morita type between $A$ and $B$, which we denote by $\Phi$.
Since $eAe$ and $fBf$
are self-injective algebras, the functor $G$ is clearly an
almost $\nu$-stable derived equivalence, and therefore
$f\Delta^0e\otimes_{eAe}-$ induces a stable equivalence of
Morita type, say $\Phi_1$, between $eAe$ and $fBf$.

For each $eAe$-module $X$, using the fact that $\Delta^0\otimes_AAe\in\add(Bf)$, we have the following isomorphisms in $\modcat{B}$:
$$Bf\otimes_{fBf}(f\Delta^0e\otimes_{eAe}X)\simeq \big(Bf\otimes_{fBf}fB\otimes_B(\Delta^0\otimes_AAe)\big)\otimes_{eAe}X\simeq \Delta^0\otimes_AAe\otimes_{eAe}X.$$
This implies that the functors $\Phi\lambda$ and $\lambda\Phi_1$ are naturally isomorphic, where the functor $\lambda$ was described in Section \ref{sect2.3}. Thus the desired commutative diagram in Lemma \ref{LemTiltStM} exists and the statements $(1)$ and $(2)$ then follow by definition.

$(3)$ Since ${}_B\Delta^i\in\add(Bf)$ for all $i>0$, the term $\Theta^{-i}=\Hom_B(\Delta^i, {}_BB)\in\add(fB)$ as a right $B$-module for all $i>0$. Now let $S$ be a simple $A$-module with $eS = 0$, that is, $eA\otimes_AS=0$. By the definition of $\cpx{\Delta}$ and $\cpx{\Theta}$, there is an isomorphism $\cpx{\Theta}\otimesP_B\cpx{\Delta}\otimesP_AS\simeq S$ in $\Db{A}$. Further, we have the following isomorphisms in $\Db{A}$:
$$\begin{array}{rl}
 S& \; \simeq \cpx{\Theta}\otimesP_B\cpx{\Delta}\otimesP_AS  \\
   &\simeq \cpx{\Theta}\otimesP_B(\Delta^0\otimes_AS) \quad (\mbox{because }\Delta^i_A\in\add(eA) \mbox{ for all }i>0)\\
   &\simeq {\Theta}^0\otimes_B\Delta^0\otimes_AS \quad (\mbox{because }\Theta^i_B\in\add(fB) \mbox{ for all }i<0 \mbox{ and } fB\otimes_B\Delta^0\in\add(eA)).\\
\end{array}$$
By the proof of \cite[Theorem 5.3]{HuXi3}, the bimodules $\Delta^0$ and $\Theta^0$ define a stable equivalence of Morita type between $A$ and $B$. Similar to the proof of Lemma \ref{LemPropStM} (5), we see that $\Phi(S)=\Delta^0\otimes_AS$ is a simple $B$-module. Morevoer, $f\cdot\Phi(S)\simeq fB\otimes_B\Delta^0\otimes_AS=0$ since $fB\otimes_B\Delta^0\in\add(eA)$ and $eA\otimes_AS=0$.

   (4) Using the two-sided tilting complex $\cpx{\Theta}=\HomP_B(\cpx{\Delta}, B)$, we can proceed the proof
   of (4) similarly as we have done in (3). $\square$

\medskip
In the following, we shall construct a Morita equivalence from a $\nu$-stable idempotent together with an arbitrary stable equivalence of Morita type induced from a derived equivalence.

\begin{Prop} Let $A$ be an algebra and $e$ be a $\nu$-stable idempotent element in $A$,
and let $\Delta$ be a self-injective algebra. Suppose that $\Xi:
\stmodcat{eAe}\ra\stmodcat{\Delta}$ is a stable equivalence of
Morita type induced by a derived equivalence. Then there is another
algebra $B$ (not necessarily isomorphic to $A$), a $\nu$-stable idempotent element $f$ in $B$, and a commutative
diagram of functors:
$$\xy
  (0,0)*+{\stmodcat{fBf}}="c",
  (25,0)*+{\stmodcat{eAe}}="d",
  (0,15)*+{\stmodcat{B}}="a",
  (25,15)*+{\stmodcat{A}}="b",
(50,0)*+{\stmodcat{\Delta}}="e",
  {\ar^{\lambda} "c";"a"},
  {\ar^{\lambda} "d";"b"},
  {\ar^{\Phi}, "a";"b"},
  {\ar^{{\Phi_1}}, "c";"d"},
    {\ar^{{\Xi}}, "d";"e"},
\endxy$$
such that

$(1)$ $\Phi$ is a stable equivalence of Morita type induced by an
iterated almost $\nu$-stable derived equivalence.

$(2)$ $\Phi_1$ is a stable equivalence of Morita type and
$\Xi\circ\Phi_1$ lifts to a Morita equivalence.

$(3)$ $\Phi(T)$ is isomorphic in $\stmodcat{A}$ to a simple
$A$-module $T'$ with $e\cdot T'=0$ for all simple $B$-modules $T$
with $f\cdot T=0$.

$(4)$ $\Phi^{-1}(S)$ is isomorphic in $\stmodcat{B}$ to a simple
$B$-module $S'$ with $f\cdot S'=0$ for all simple $A$-modules $S$
with $e\cdot S=0$. \label{PropExtStM}
\end{Prop}

{\it Proof.} Suppose that the stable equivalence $\Xi$ is induced by a standard derived equivalence
$F: \Db{eAe}\ra \Db{\Delta}$. Then there is an integer $m\le 0$
such that $[m]\circ F$ is an almost $\nu$-stable derived equivalence. Observe that the shift functor $[1]:
\Db{\Delta}\ra\Db{\Delta}$ is isomorphic to the standard derived
equivalence $(\Delta[1])\otimesL_{\Delta}-$. Thus, the derived
equivalence $[m]\circ F$ is standard, and
consequently $\Xi$ can be written as a composite $\Xi\simeq \Xi_2\circ\Xi_1$ of stable equivalences $\Xi_1$ and $\Xi_2$ of Morita type such that
$\Xi_1$ is induced by $[m]\circ F:\Db{eAe}\ra\Db{\Delta}$ and $\Xi_2$
is induced by $[-m]: \Db{\Delta}\ra\Db{\Delta}$.

Let $\cpx{X}$ be a tilting complex over $eAe$ associated to $[m]\circ F$. Then
$X^i=0$ for all $i>0$. Set $\cpx{Q}:=Ae\otimesP_{eAe}\cpx{X}$.
Then $\cpx{Q}$ satisfies all conditions in Lemma \ref{LemTiltStM} since $e\cpx{Q}\simeq \cpx{X}$ is a tilting complex over $eAe$.
Hence, by Lemma \ref{LemTiltStM}, there is an algebra $B'$ and a
$\nu$-stable idempotent element $f'$ in $B'$, together with a commutative diagram:
  $$\xy
    (0,30)*+{\stmodcat{B'}}="b",
    (0,15)*+{\stmodcat{f'B'f'}}="fbf",
    (0.1,0.1)*+{\Db {f'B'f'}}="DbfBf",
    (30,30)*+{\stmodcat{A}}="a",
    (30,15)*+{\stmodcat{eAe}}="eae",
    (30, 0)*+{\Db{eAe}}="Dbeae",
    (60, 15)*+{\stmodcat{\Delta}}="lam",
    (60, 0)*+{\Db{\Delta}}="Dblam",
    {\ar^{[m]\circ F} "Dbeae"; "Dblam"},
    {\ar_{G_1}  "Dbeae"; "DbfBf"},
    {\ar^{\eta_{\Delta}} "Dblam"; "lam"},
    {\ar^{\eta_{eAe}} "Dbeae"; "eae"},
    {\ar^{\eta_{f'B'f'}} "DbfBf"; "fbf"},
    {\ar^{\Xi_1} "eae";"lam"},
    {\ar_(.4){\Phi_1'} "eae"; "fbf"},
    {\ar_{\Phi'} "a"; "b"},
    {\ar^{\lambda} "eae"; "a"},
    {\ar^{\lambda} "fbf"; "b"},
  \endxy$$
such that $\Phi'$ is a stable equivalence of Morita type induced by
a standard, almost $\nu$-stable derived equivalence, and that $G_1$
is a standard derived equivalence with $\cpx{X}$ as an associated
tilting complex. Thus $f'B'f'$ is a tilting complex associated to
the derived equivalence $[m]\circ F\circ G_1^{-1}$. This means that
$(f'B'f')[m]$ is a tilting complex associated to $F\circ G_1^{-1}$.
By the dual version of Lemma \ref{LemTiltStM}, there is an algebra
$B$ and a $\nu$-stable idempotent element $f$ in $B$, together with
a commutative diagram
 $$\xy
    (0,30)*+{\stmodcat{B}}="b",
    (0,15)*+{\stmodcat{fBf}}="fbf",
    (0,0)*+{\Db{fBf}}="DbfBf",
    (30, 30)*+{\stmodcat{B'}}="a",
    (30,15)*+{\stmodcat{f'B'f'}}="eae",
    (30,0)*+{\Db{f'B'f'}}="Dbeae",
    (70, 15)*+{\stmodcat{\Delta}}="lam",
    (70, 0)*+{\Db{\Delta}}="Dblam",
    {\ar^{F\circ G_1^{-1}} "Dbeae"; "Dblam"},
    {\ar_{G_2}  "Dbeae";"DbfBf"},
    {\ar^{\eta_{\Delta}} "Dblam"; "lam"},
    {\ar^{\eta_{f'B'f'}} "Dbeae"; "eae"},
    {\ar^{\eta_{fBf}} "DbfBf"; "fbf"},
    {\ar^(.55){\Xi_2\Xi_1{(\Phi_1')}^{-1}} "eae";"lam"},
    {\ar_{\Phi''_1} "eae"; "fbf"},
    {\ar_{\Phi''} "a"; "b"},
    {\ar^{\lambda} "eae"; "a"},
    {\ar^{\lambda} "fbf";"b"},
  \endxy$$
 such that $\Phi''$ is a stable equivalence of Morita type with
 $(\Phi'')^{-1}$ induced by a standard almost $\nu$-stable derived
 equivalence, and that $G_2$ is a standard derived equivalence with
 $(f'B'f')[m]$ as an associated tilting complex.

 Now we define $\Phi:=(\Phi')^{-1}\circ(\Phi'')^{-1}$,
 $\Phi_1:={(\Phi'_1)}^{-1}\circ{(\Phi''_1)}^{-1}$, we get the following commutative
 diagram
  $$\xy
    (0, 30)*+{\stmodcat{B}}="b",
    (0,15)*+{\stmodcat{fBf}}="fbf",
    (0,0)*+{\Db{fBf}}="DbfBf",
    (30, 30)*+{\stmodcat{A}}="a",
    (30,15)*+{\stmodcat{eAe}}="eae",
    (30, 0)*+{\Db{eAe}}="Dbeae",
    (70, 15)*+{\stmodcat{\Delta}}="lam",
    (70, 0)*+{\Db{\Delta}}="Dblam",
    {\ar^{F} "Dbeae"; "Dblam"},
    {\ar^{G_1^{-1}\circ G_2^{-1}} "DbfBf"; "Dbeae"},
    {\ar^{\eta_{\Delta}} "Dblam"; "lam"},
    {\ar^{\eta_{eAe}} "Dbeae"; "eae"},
    {\ar^{\eta_{fBf}} "DbfBf"; "fbf"},
    {\ar^{\Xi_2\circ\Xi_1} "eae";"lam"},
    {\ar^{\Phi_1}  "fbf"; "eae"},
    {\ar^{\Phi} "b"; "a"},
    {\ar^{\lambda} "eae"; "a"},
    {\ar^{\lambda} "fbf";"b"},
  \endxy$$
By the above discussion, we see that $fBf$ is a tilting complex
associated to $F\circ G_1^{-1}\circ G_2^{-1}$. Hence the derived equivalence $F\circ
G_1^{-1}\circ G_2^{-1}$ is induced by a Morita equivalence. Consequently, the stable equivalence
$\Xi\circ\Phi_1\simeq\Xi_2\circ\Xi_1\circ\Phi_1$ lifts
to a Morita equivalence. Thus (1) and (2) follow.  Now, (3) and (4) follow
easily from the above diagram and Lemma \ref{LemTiltStM} (3)-(4).
$\square$

\subsection{Proof of Theorem \ref{ThmLift}\label{sectionFromStMToDer}}

With the above preparations, we now give a proof of Theorem \ref{ThmLift}.

\medskip
{\it\parindent=0pt Proof of Theorem \ref{ThmLift}.} By Lemma \ref{LemRes}, there is a stable equivalence
$\Phi_1:\stmodcat{eAe}\ra\stmodcat{fBf}$ of Morita type such that the following diagram of functors
$$\xy
  (0,0)*+{\stmodcat{eAe}}="c",
  (25,0)*+{\stmodcat{fBf}}="d",
  (0,15)*+{\stmodcat{A}}="a",
  (25,15)*+{\stmodcat{B}}="b",
  {\ar^{\lambda} "c";"a"},
  {\ar^{\lambda} "d";"b"},
  {\ar^{\Phi}, "a";"b"},
  {\ar^{\Phi_1}, "c";"d"},
\endxy$$
is commutative up to isomorphism. Note that the functor
$\Phi_1$ is uniquely determined up to isomorphism.

By Proposition \ref{PropExtStM}, we can find an algebra $B'$, a
$\nu$-stable idempotent element $f'$ in $B'$ and stable equivalences $\Phi':
\stmodcat{B'}\ra \stmodcat{A}$ and $\Phi_1': \stmodcat{f'B'f'}\ra
  \stmodcat{eAe}$ together with a commutative diagram
  $$\xy
    (0,0)*+{\stmodcat{f'B'f'}}="fbfl",
    (0,15)*+{\stmodcat{B'}}="bl",
    (30,0)*+{\stmodcat{eAe}}="eae",
    (30,15)*+{\stmodcat{A}}="a",
    (60,0)*+{\stmodcat{fBf}}="fbfr",
    (60,15)*+{\stmodcat{B}}="br",
    {\ar^{\Phi'} "bl";"a"},
    {\ar^{\Phi_1'} "fbfl";"eae"},
    {\ar^{\Phi} "a";"br"},
    {\ar^{\Phi_1} "eae";"fbfr"},
    {\ar^{\lambda} "fbfl";"bl"},
    {\ar^{\lambda} "eae";"a"},
    {\ar^{\lambda} "fbfr";"br"},
  \endxy$$
such that $\Phi_1\circ\Phi_1'$ lifts to a Morita equivalence and
$\Phi'$ is induced by an iterated almost $\nu$-stable derived
equivalence. Moreover, for all simple $B'$-modules $S'$ with $f'\cdot
S'=0$, the module $\Phi'(S')$ is isomorphic to a simple module
$S$ with $e\cdot S=0$, and dually, for all simple
$A$-modules $V$ with $e\cdot V=0$,  the module ${\Phi'}^{-1}(V)$ is isomorphic in $\stmodcat{B'}$ to a
simple $B'$-module $V'$ with $f'\cdot V'=0$. By the assumptions (1) and (2) in
Theorem \ref{ThmLift}, the $B$-module $(\Phi\circ\Phi')(S')$ is
isomorphic in $\stmodcat{B}$ to a simple module for each simple
$B'$-module $S'$ with $f'\cdot S'=0$; and the $B'$-module
${\Phi'}^{-1}\circ\Phi^{-1}(V)$ is isomorphic to a simple
$B'$-module for each simple $B$-module $V$ with $f\cdot V=0$. It
follows from Proposition \ref{PropLiftMorita} that $\Phi\circ\Phi'$
lifts to a Morita equivalence, and consequently the functor $\Phi$
is induced by an iterated almost $\nu$-stable equivalence since
$\Phi'$ is induced by an iterated almost $\nu$-stable derived
equivalence. This finishes the proof of Theorem \ref{ThmLift}. $\square$

\medskip
Let us remark that every stable equivalence of Morita type between
two algebras $A$ and $B$ can be ``restricted" to a stable
equivalence of Morita type between $eAe$ and $fBf$ for some
$\nu$-stable idempotent elements $e\in A$ and $f\in B$. There are two typical
ways to implement this point:

(i) For each algebra $A$, there is an associated self-injective algebra (see Subsection 2.3), which we denote by $\Delta_A$.
The result \cite[Theorem 4.2]{DugasVilla} shows that if $A/\rad(A)$ and $B/\rad(B)$ are separable then every stable
equivalence of Morita type between $A$ and $B$ restricts to a stable
equivalence of Morita type between the associated self-injective
algebras $\Delta_A$ and $\Delta_B$.

(ii) Under the setting of Lemma \ref{LemPropStM}, let $e$ be an
idempotent element of $A$ such that $\add(Ae)=\add({}_AP)$. Let $f$ be
defined similarly such that $\add(Bf)=\add({}_BQ)$. Then it follows from
Lemma \ref{LemPropStM} (2) that the idempotent elements $e$ and $f$ are
$\nu$-stable. By Lemma \ref{LemPropStM} (3) and \cite[Theorem
 1.2]{ChenPanXi}, the given stable equivalence of Morita type between $A$
and $B$ in Lemma \ref{LemPropStM} restricts to a stable equivalence
of Morita type between $eAe$ and $fBf$.

\medskip
As an immediate consequence of Theorem \ref{ThmLift}, we have the following corollary.

\begin{Koro} \label{CorASelfinj}
Let $A$ and $B$ be algebras without
semisimple direct summands such that $A/\rad(A)$ and $B/\rad(B)$ are
separable. Suppose that $\Phi$ is a stable equivalence of Morita
type between $A$ and $B$, and that $\Phi_1$ is the restricted stable
equivalence of $\Phi$ between the associated self-injective algebras
$\Delta_A$ and $\Delta_B$. If $\Phi_1$ lifts to a derived
equivalence between $\Delta_A$ and $\Delta_B$, then $\Phi$ lifts to an iterated almost $\nu$-stable
derived equivalence between $A$ and $B$.
\end{Koro}

{\it Proof.} By definition,  $\Delta_A=eAe$ for some idempotent $e$ in $A$ with $\add(Ae)=\stp{A}$, and $\Delta_B=fBf$ for some idempotent $f$ in $B$ with $\add(Bf)=\stp{B}$. Suppose that $_AM_B$ and ${}_BN_A$ are two bimodules without projective direct summands and define a stable equivalence of Morita type between $A$ and $B$ such that $\Phi$ is induced by $_BN\otimes_A-$.  Assume that ${}_AM\otimes_BN_A\simeq A\oplus {}_AP_A$ and ${}_BN\otimes_AM_B\simeq B\oplus {}_BQ_B$ as bimodules.

We first show that $N\otimes_AAe\in\add({}_BBf)$ and $M\otimes_BBf\in\add({}_AAe)$. By the proof of Lemma \ref{LemPropStM} (2), we have $\nu^i_A(N\otimes_A Ae)\simeq N\otimes_A(\nu_A^i(Ae))$ for all $i\geq 0$. Note that $\nu_A^i(Ae)$ is projective for all $i\geq 0$ since $Ae\in \stp{A}$. This implies that $\nu^i_B(N\otimes_A Ae)$ is projective for all $i\geq 0$, that is, $N\otimes_AAe\in\stp{B}=\add({}_BBf)$. Similarly, we have $M\otimes_BBf\in\add({}_AAe)$.

Let $S$ be a simple $A$-module with $e\cdot S=0$. By Lemma
\ref{LemPropStM} (2), the module ${}_AP$ is in $\stp{A}$, which is
exactly $\add({}_AAe)$. Hence $\Hom_A(P, S)=0$ and consequently
$\Phi(S)=N\otimes_AS$ is a simple $B$-module by Lemma
\ref{LemPropStM} (5). Moreover,
$$f\cdot \Phi(S)=\Hom_B(Bf, N\otimes_AS)\simeq \Hom_A(M\otimes_BBf, S)=0$$
since $M\otimes_BBf\in\add(Ae)$. Similarly, for each simple
$B$-module $V$ with $f\cdot V=0$, the $A$-module $\Phi^{-1}(V)$ is
simple with $e\cdot\Phi^{-1}(V)=0$. Now, the corollary follows from
Theorem \ref{ThmLift}. $\square$

\medskip
In the next section we will find out a class of algebras for which
$\Phi_1$ can be lifted to a derived equivalence.

\section{Frobenius-finite algebras: Proof of Theorem \ref{ThmRepFin}\label{sect4.3}}

Corollary \ref{CorASelfinj} shows that the associated self-injective
algebra of a given algebra may be of prominent importance in lifting
stable equivalences of Morita type to derived equivalences. Based on
this point of view, we shall introduce, in this section, a class of algebras, called
Frobenius-finite algebras, and discuss their basic properties. With these preparations in hand, we then prove
Theorem \ref{ThmRepFin}.

\subsection{Frobenius-finite algebras \label{sect5.1}}

\begin{Def}
A finite-dimensional $k$-algebra is said to be {\bf Frobenius-finite} if
its associated self-injective algebra is representation-finite, and
Frobenius-free if its associated self-injective algebra is zero.
\end{Def}
Clearly, Frobenius-free algebras are Frobenius-finite, and
representation-finite algebras are Frobenius-finite. Moreover, the ubiquity of Frobenius-finite algebras is guaranteed by the next propositions.

Before we present methods to product Frobenius-finite algebras, let us recall the definition of Auslander-Yoneda algebras introduced in \cite{HuXi4}. A subset $\Theta$ of $\mathbb{N}$ is called an \emph{admissible subset} if $0\in\Theta$ and if, for each $l,m,n\in\Theta$ with $l+m+n\in\Theta$, we have $l+m\in\Theta$ if and only if $m+n\in\Theta$. There are many admissible subsets of $\mathbb{N}$. For example, for each $n\in\mathbb{N}$, the subsets $\{xn\mid x\in\mathbb{N}\}$ and $\{0,1,2,\cdots,n\}$ of $\mathbb{N}$ are admissible.

Let $\Theta$ be an admissible subset of $\mathbb{N}$, and let ${\cal T}$ be a triangulated $k$-category. There is a bifunctor
$$\E{\cal T}{\Theta}{-}{-}: {\cal T}\times {\cal T}\lra \Modcat{k}$$
$$(X, Y)\mapsto \E{\cal T}{\Theta}{X}{Y}:=\bigoplus_{i\in\Theta}\Hom_{\cal T}(X, Y[i])$$
with composition given in an obvious way (for details, see \cite[Subsection 3.1]{HuXi4}). In particular, if $f\in\Hom_{\cal T}(X, Y[i])$ and $g\in\Hom_{\cal T}(Y, Z[j])$, then the composite $f\cdot g=f(g[i])$ if $i+j\in\Theta$, and $f\cdot g=0$ otherwise. In this way, for each object $M\in{\cal T}$, we get an associated algebra $\E{\cal T}{\Theta}{M}{M}$, which is simply denoted by $\Ex{\cal T}{\Theta}{M}$ and called $\Theta$-\emph{Auslander-Yoneda algebra} of $M$. If ${\cal T}=\Db{A}$ for some algebra $A$, we denote $\E{\Db{A}}{\Theta}{X}{Y}$ by $\E{A}{\Theta}{X}{Y}$, and $\Ex{\Db{A}}{\Theta}{M}$ by $\Ex{A}{\Theta}{M}$ for all $X, Y, M\in\Db{A}$.

The following proposition shows that Frobenius-finite algebras can be constructed from
generator-cogenerators. Thus there are plenty of Frobenius-finite algebras. Recall that an $A$-module $M$ is called a \emph{generator} in $\modcat{A}$ if $\add(M)$ contains ${}_AA$;
a {\em generator-cogenerator} in $\modcat{A}$ if $\add(M)$ contains both
${}_AA$ and ${}_AD(A)$; and a \emph{torsionless} module if it is a submodule of a projective module.

\begin{Prop}
$(1)$ Let $M$ be a generator-cogenerator over a Frobenius-finite algebra $A$. Then $\End_A(M)$ is Frobenius-finite.
In particular, Auslander algebras are Frobenius-finite.

$(2)$ Let $M$ be a torsionless generator over a Frobenius-finite algebra $A$. Suppose that $\Theta$ is a finite admissible subset of $\mathbb{N}$ and that $\Ext_A^i(M, A)=0$ for all $0\neq i\in\Theta$. Then $\Ex{A}{\Theta}{M}$ is Frobenius-finite. In particular, if $A$ is a representation-finite self-injective algebra, then $\Ex{A}{\Theta}{A\oplus X}$ is Frobenius-finite for each $A$-module $X$ and for arbitrary finite admissible subset $\Theta$ of $\mathbb{N}$.

$(3)$ If $A$ and $B$ are Frobenius-finite algebras and $_BM_A$ is a
bimodule, then the triangular matrix algebra
$\left[\begin{smallmatrix}A & 0\\M & B\end{smallmatrix}\right]$ is
Frobenius-finite. More generally, if $A_1, \cdots, A_m$ are a family of Frobenius-finite algebras and if $M_{ij}$ is an $A_i$-$A_j$-bimodule for all $1\le j<i\le m$, then the
triangular matrix algebra of the form
$$\left[\begin{array}{cccc}A_1 \\M_{21} & A_2\\
      \vdots & \vdots & \ddots \\
      M_{m1} & M_{m2} &\cdots & A_m\end{array}\right]$$
is Frobenius-finite.

$(4)$ If $A=A_0\oplus A_1\oplus \cdots\oplus A_n$ is
an $\mathbb N$-graded algebra with $A_0$ Frobenius-finite, then the Beilinson-Green algebra
$$\Lambda_m:=\left[\begin{array}{cccc}A_0 \\ A_1 & A_0\\
      \vdots & \ddots & \ddots \\
      A_m & \cdots & A_1  & A_0\end{array}\right]$$
is Frobenius-finite for all $1\le m\le n$.

\end{Prop}

{\it Remark.} The triangular matrix algebra of a graded algebra $A$ in (4) seems first to appear
in a paper by Edward L. Green in 1975. A special case of this
kind of algebras appeared in a paper by A. A. Beilinson in 1978, where he described the derived category of coherent sheaves over $\mathbb{P}^n$ as the one of this triangular matrix algebra.
Perhaps it is more appropriate to name this triangular matrix algebra as the
\emph{Beilinson-Green} algebra of $A$.

\medskip
{\it Proof.} (1) We set $\Lambda:=\End_A(M)$. Since $M$ is a generator-cogenerator for $A$-mod, every indecomposable projective-injective
$\Lambda$-module is of the form $\Hom_A(M, I)$ with $I$ an
indecomposable injective $A$-module. Moreover, for each
projective $A$-module $P$, there is a natural isomorphism
$\nu_{\Lambda}\Hom_A(M, P)\simeq \Hom_A(M, \nu_AP)$.
This implies that $\Hom_A(M, P)\in\stp{\Lambda}$ for all
$P\in\stp{A}$. Now let $I$ be an indecomposable injective $A$-module
such that $\Hom_A(M, I)$ lies in $\stp{\Lambda}$. Then it follows
from $\nu_{\Lambda}^{-1}\Hom_A(M, I)\simeq \Hom_A(M, \nu_A^{-1}I)$
that $\Hom_A(M,\nu^{-1}_AI)$ lies in $\stp{\Lambda}$. Consequently,
the $\Lambda$-module $\Hom_A(M,\nu^{-1}_AI)$ is injective, and
therefore the $A$-module $\nu_A^{-1}I$ is projective-injective.
Applying $\nu_{\Lambda}^{-1}$ repeatedly, one sees that $\nu_A^{i}I$
is projective-injective for all $i<0$. This implies $I\in\stp{A}$.
Hence the restriction of the functor $\Hom_A(M, -): \add(_AM)\ra
\pmodcat{\Lambda}$ to the category $\stp{A}$ is an equivalence
between $\stp{A}$ and $\stp{\Lambda}$. Consequently, the associated
self-injective algebras $\Delta_A$ and $\Delta_{\Lambda}$ are Morita
equivalent. Thus (1) follows.

\medskip
(2) For convenience, we set
$\Lambda:=\Ex{A}{\Theta}{M}=\bigoplus_{i\in\Theta}\Lambda_i$, where
$\Lambda_i:=\Hom_{\Db{A}}(M, M[i])$. We also identify $\Ext_A^i(U,
V)$ with $\Hom_{\Db{A}}(U, V[i])$ for all $A$-modules $U, V$ and all
integers $i$. Observe that $\rad(\Lambda)=\rad(\Lambda_0)\oplus
\Lambda_{+}$, where $\Lambda_{+}:=\bigoplus_{0\neq
i\in\Theta}\Lambda_i$.

We  shall prove that $\stp{A}$ and $\stp{\Lambda}$ are equivalent.
Let $Y$ be an indecomposable non-projective direct summand of $M$.
We claim that $\E{A}{\Theta}{M}{Y}$ cannot be in $\stp{\Lambda}$.
Suppose contrarily that this is false and
$\E{A}{\Theta}{M}{Y}\in\stp{\Lambda}$. Then the $\Lambda$-module
$\E{A}{\Theta}{M}{Y}$ must be indecomposable projective-injective.
Now, we have to consider the following two cases:

(a) $\bigoplus_{0\neq i\in\Theta}\Ext_A^i(M, Y)=0$. Since $Y$ is
torsionless, there is an injective $A$-module homomorphism $f: Y\ra
A^n$. Applying $\Hom_A(M, -)$ results in an injective map $\Hom_A(M,
f): \Hom_A(M, Y)\ra \Hom_A(M, A^n)$. Together with the assumption
that $\Ext_A^i(M, A)=0$ for all $0\neq i\in\Theta$, we see that
$\E{A}{\Theta}{M}{Y}=\Hom_A(M, Y)$, $\E{A}{\Theta}{M}{A^n}=\Hom_A(M,
A^n)$ and $\E{A}{\Theta}{M}{f}=\Hom_A(M, f)$. This implies that
$\E{A}{\Theta}{M}{f}: \E{A}{\Theta}{M}{Y}\ra\E{A}{\Theta}{M}{A^n}$
is an injective map and must splits. It follows that $Y$ must be a
direct summand of $A^n$. This is a contradiction.

(b)  $\bigoplus_{0\neq i\in\Theta}\Ext_A^i(M, Y)\neq 0$. Let $m\neq
0$ be the maximal integer in $\Theta$ with $\Ext_A^m(M, Y)\neq 0$.
Then $\Lambda_{+}\Ext_A^m(M, Y)=0$, and consequently
$\rad(\Lambda)\soc_{\Lambda_0}\big(\Ext_A^m(M, Y)\big)=0$. Hence
$\soc_{\Lambda_0}\big(\Ext_A^m(M, Y)\big)=\Lambda\cdot
\soc_{\Lambda_0}\big(\Ext_A^m(M, Y)\big)$ is a $\Lambda$-submodule
of $\soc_{\Lambda}\big(\E{A}{\Theta}{M}{Y}\big)$. Next, we show that
$\soc_{\Lambda_0}\big(\Hom_A(M, Y)\big)$ is also a
$\Lambda$-submodule of
$\soc_{\Lambda}\big(\E{A}{\Theta}{M}{Y}\big)$. Let $g: M\ra Y$ be in
$\soc_{\Lambda_0}\big(\Hom_A(M, Y)\big)$. Suppose that $M=M_p\oplus
X$ such that $M_p$ is projective and $X$ contains no projective
direct summands. Now for each $x\in X$, there are indecomposable
projective modules $P_j$, $1\le j\le s$ and homomorphisms $h_j:
P_j\ra X$, which must be radical maps,  such that $x=\sum_{j=1}^s
(p_j)h_j$ for some $p_j\in P_j$ with $j=1, \cdots, s$. Since $M$ is
a generator over $A$, the module $P_j$ is isomorphic to a direct
summand of $M$. Thus, we get a map $\tilde{h}_j: M\ra P_j\raf{h_j}
X\hookrightarrow M$, which is in $\rad(\Lambda_0)$ for all $j$, and
the composite $\tilde{h_j}g$ has to be zero. This implies that the
image of $x$ under $g$ is $0$, and consequently $g|_{X}=0$.
%there is an indecomposable projective $A$-module $P$ together with an $A$-module homomorphism $h: P\ra X$, which must be a radical map, such that $x$ lies in the image of $h$. Since $M$ is a generator over $A$, $P$ is isomorphic to a direct summand of $M$. Thus, we get a map $\tilde{h}: M\ra P\raf{h} X\hookrightarrow M$, which is in $\rad(\Lambda_0)$, and the composite $\tilde{h}g$ has to be zero. This implies that the image of $x$ under $g$ is $0$, and consequently $g|_{X}=0$.
Let $\pi: M\ra M_p$ be the canonical projection. Then we have $g=\pi
g'$ for some $g': M_p\ra Y$. For each $t: M\ra M[i]$  in $\Db{A}$
with $0\neq i\in\Theta$, the composite $t\cdot g=t
(g[i])=t(\pi[i])(g'[i])$. Since $\Ext_A^i(M, A)=0$, we have
$\Ext_A^i(M, M_p)=0$, and consequently $t(\pi[i])=0$. Hence $t\cdot
g=0$, and therefore $\Lambda_{+}\cdot \soc_{\Lambda_0}\big(\Hom_A(M,
Y)\big)=0$. It follows that $\rad(\Lambda)\cdot
\soc_{\Lambda_0}\big(\Hom_A(M, Y)\big)=0$ and that
$\soc_{\Lambda_0}\big(\Hom_A(M, Y)\big)=\Lambda\cdot
\soc_{\Lambda_0}\big(\Hom_A(M, Y)\big)$ is a $\Lambda$-submodule of
$\soc_{\Lambda}\big(\E{A}{\Theta}{M}{Y}\big)$. Thus, we have shown
that the $\Lambda$-module $\soc_{\Lambda_0}\big(\Hom_A(M,
Y)\big)\oplus \soc_{\Lambda_0}\big(\Ext_A^m(M, Y)\big)$ is contained
in $\soc_{\Lambda}\big(\E{A}{\Theta}{M}{Y}\big)$. This shows that
$\soc_{\Lambda}\big(\E{A}{\Theta}{M}{Y}\big)$ cannot be simple and
$\E{A}{\Theta}{M}{Y}$ cannot be indecomposable injective. This is
again a contradiction.

Thus, we have shown that every indecomposable projective $\Lambda$-module in $\stp{\Lambda}$ has to be of the form $\E{A}{\Theta}{M}{P}$ for some indecomposable projective $A$-module $P$. Suppose $\E{A}{\Theta}{M}{P}\in \stp{\Lambda}$. We shall prove $P\in\stp{A}$. In fact, by \cite[Lemma 3.5]{HuXi4}, we have $\nu_{\Lambda}\E{A}{\Theta}{M}{P}\simeq \E{A}{\Theta}{M}{\nu_AP}$. It follows from definition that $\nu_{\Lambda}\E{A}{\Theta}{M}{P}$ is again in $\stp{\Lambda}$. This means that there is an isomorphism $\E{A}{\Theta}{M}{\nu_AP}\simeq \E{A}{\Theta}{M}{P'}$ for some indecomposable projective $A$-module $P'$. Since $\Ext_A^i(M, A)=0$ for all $0\neq i\in\Theta$ and since $\nu_AP$ is injective, we have $\Hom_A(M, \nu_AP)=\E{A}{\Theta}{M}{\nu_AP}\simeq \E{A}{\Theta}{M}{P'}=\Hom_A(M, P')$. Hence $\nu_AP\simeq P'$ is projective by Lemma \ref{lemEndMproj}, Repeatedly, we see that $\nu_A^iP$ is projective for all $i>0$, that is, $P\in\stp{A}$. Conversely, let $P$ be an indecomposable module in $\stp{A}$. Due to the isomorphism $\nu_{\Lambda}\E{A}{\Theta}{M}{P}\simeq \E{A}{\Theta}{M}{\nu_AP}$, the $\Lambda$-module $\E{A}{\Theta}{M}{P}$ belongs to $\stp{\Lambda}$. Thus, the functor $\E{A}{\Theta}{M}{-}$ induces an equivalence from $\stp{A}$ to $\stp{\Lambda}$. Hence the associated self-injective algebras $\Delta_A$ and $\Delta_{\Lambda}$ are Morita equivalent, and (2) follows.

(3) Set $\Lambda:=\left[\begin{smallmatrix}A & 0\\M &
B\end{smallmatrix}\right]$. Then each $\Lambda$-module can be interpreted
as a triple $({}_AX, {}_BY, f)$ with $X\in\modcat{A}$,
$Y\in\modcat{B}$ and $f: {}_BM\otimes_AX\ra {}_BY$ a $B$-module
homomorphism. Let $({}_AX, {}_BY, f)$ be an indecomposable
$\Lambda$-module in $\stp{\Lambda}$. Then $({}_AX, {}_BY, f)$ is
projective-injective with $\nu_{\Lambda}({}_AX, {}_BY,
f)\in\stp{\Lambda}$. By \cite[p.76, Proposition 2.5]{AusReiten},
there are two possibilities:

\smallskip
(i) ${}_BY=0$ and ${}_AX$ is an indecomposable projective-injective
$A$-module with $M\otimes_AX=0$;

(ii) ${}_AX=0$ and ${}_BY$ is an indecomposable projective-injective
$B$-module with $\Hom_B(M, Y)=0$.

\smallskip
{\parindent=0pt Now} we assume (i). Then $\nu_{\Lambda}(X, 0,
0)\simeq (\nu_AX, 0, 0)$ is still in $\stp{\Lambda}$. This implies that $\nu_A^iX$ is projective-injective for all
$i\geq 0$, and therefore $X\in\stp{A}$. Similarly, if we assume (ii), then
$Y\in\stp{B}$.  Thus, we can assume that $\{(X_1, 0, 0), \cdots,
(X_r, 0, 0), (0, Y_1, 0), \cdots, (0, Y_s, 0)\}$ is a complete set
of non-isomorphic indecomposable modules in $\stp{\Lambda}$ with both
$X_i\in\stp{A}$ and $Y_j\in\stp{B}$ for all $i$ and $j$. Then the
associated self-injective algebra
$$\Delta_{\Lambda}:=\End_{\Lambda}\big(\bigoplus_{i=1}^r(X_i, 0,
0)\oplus \bigoplus_{i=1}^s(0, Y_i, 0)\big)\simeq
\End_A(\bigoplus_{i=1}^rX_i)\times \End_B(\bigoplus_{i=1}^sY_i)$$ is
representation-finite since both $A$ and $B$ are Frobenius-finite.

(4) This is an immediate consequence of (3). $\square$

\medskip
In the following, we shall show that Frobenius-finite algebras can be obtained by Frobenius extensions.

Suppose that $B$ is a subalgebra of an algebra $A$. We denote by $F$ the induction functor ${}_AA\otimes_B-: \modcat{B}\ra \modcat{A}$ and by $H$ the restriction functor ${}_B(-): \modcat{A}\ra\modcat{B}$. Observe that for any $k$-algebra $C$, the functor $F$ is also a functor from $B$-$C$-bimodules to $A$-$C$-bimodules and $H$ is also a functor from $A$-$C$-bimodules to $B$-$C$-bimodules.

\begin{Prop}
Let $B$ be a subalgebra of an algebra $A$. Suppose that the extension $B\hookrightarrow A$ is Frobenius, that is, $\Hom_B(_BA,-)\simeq F$ as functors from $\modcat{B}$ to $\modcat{A}$.

$(1)$ Suppose that the extension $B\hookrightarrow A$ spits, that is, the inclusion map $B\ra A$ is a split monomorphism of $B$-$B$-bimodules.
If $A$ is Frobenius-finite, then so is $B$.

$(2)$ Suppose that the extension  $B\hookrightarrow A$ is separable, thst is, the multiplication map $A\otimes_BA\ra A$ is a split epimorphism of $A$-$A$-bimodules. If $B$ is Frobenius-finite, then so is $A$.
\label{PropalgExtension}
\end{Prop}

{\it Proof.} Clearly, $(F, H)$ is an adjoint pair. Note that, for a Frobenius extension, ${}_BA$ is a finitely generated projective module and $\Hom_B(_BA, B)\simeq A$ as $A$-$B$-bimodules (see \cite[40.21, p.423]{bw}). We first show that both $F$ and $H$ commutes with the Nakayama functors. In fact, for each $B$-module $X$, we have the following natural isomorphisms of $A$-modules:
$$\begin{array}{rl}
    \nu_{A}(F(X)) & = D\Hom_A({}_AA\otimes_BX, {}_AA_A) \\
      & \simeq D\Hom_B(X, {}_BA_A)  \quad ((F, H) \mbox{ is an adjoint pair}) \\
      & \simeq D\Hom_B(X, {}_BB\otimes_BA_A) \\
      & \simeq D\big(\Hom_B(X, B)\otimes_BA_A \big) \quad ({}_BA \mbox{ is projective})\\
      & \simeq \Hom_B(_BA_A, {}_BD(X^*)\\
      & \simeq \Hom_B\big(_BA\otimes_AA_A, {}_BD(X^*)\big)\\
      & \simeq \Hom_A\big({}_AA_A, \Hom_B(_BA_A, {}_BD(X^*)\big)\\
      & \simeq \Hom_A(_AA_A, F(\nu_B(X))) \quad (\mbox{Frobenius extension})\\
      & \simeq F(\nu_B(X)).
  \end{array}
$$
For each $A$-module $Y$, we have the following natural isomorphisms
of $B$-modules:
$$\begin{array}{rl}
    \nu_B(H(Y)) & =D\Hom_B(_BA\otimes_AY, {}_BB_B) \\
               & \simeq D\Hom_A(Y, \Hom_B(_BA, {}_BB_B))\\
     &  \simeq D\Hom_A(Y, {}_AA_B) \quad (\mbox{Frobenius extension}) \\
      & = H(\nu_A(Y)).  \\
  \end{array}
$$

Note that the functor $F$ takes projective $B$-modules to projective $A$-modules. For each projective $B$-module $P$ in $\stp{B}$, we have $\nu_A^iF(P)\simeq F(\nu_B^iP)$ is projective for all $i\geq 0$, that is, $F(P)\in\stp{A}$. Since ${}_BA$ is projective, the functor $H$ takes projective $A$-modules to projective $B$-modules. Similarly, we can show that $H(Q)$ belongs to $\stp{B}$ for all $Q\in\stp{A}$.

\smallskip
Let $e$ and $f$ be idempotents in $A$ and $B$, respectively, such that $\add(Ae)=\stp{A}$ and $\add(Bf)=\stp{B}$. Then $eAe$ and $fBf$ are the Frobenius parts of $A$ and $B$, respectively.

Note that there is an equivalence between $\modcat{fBf}$ and the full subcategory, denoted by $\mbox{mod}(Bf)$,  of $\modcat{B}$ consisting of $B$-modules $X$ that admit a projective presentation $P_1\ra P_0\ra X\ra 0$ with $P_i\in\add(Bf)$ for $i=0,1$. Similarly, the module category $\modcat{eAe}$ is equivalent to the full subcategory $\mbox{mod}(Ae)$ of $\modcat{A}$. Now for each $B$-module $X$ in $\mbox{mod}(Bf)$, let $P_1\ra P_0\ra X\ra 0$ be a presentation of $X$ with $P_0, P_1\in\add(Bf)=\stp{B}$. Applying the induction functor $F$ which is right exact, we get an exact sequence
$F(P_1)\ra F(P_0)\ra F(X)\ra 0$ with $F(P_i)$ in $\stp{A}=\add(Ae)$. This shows that $F(X)$ is in $\mbox{mod}(Ae)$ for all $X\in \mbox{mod}(Bf)$. Since the restriction of scalars functor $H$ is exact, we can deduce that $H(Y)$ lies in $\mbox{mod}(Bf)$ for all $A$-modules $Y$ in $\mbox{mod}(Ae)$.

(1) Now for each $B$-module $X$ in $\mbox{mod}(Bf)$, the assumption (1) implies that $X$ is a direct summand of $HF(X)$. If $X$ is indecomposable, then $X$ is a direct summand of $H(Y)$ for some indecomposable direct summand $Y$ of $F(X)$, which is in $\mbox{mod}(Ae)$. Thus, if $eAe$ is representation-finite, then $\mbox{mod}(Ae)$ has finitely many isomorphism classes of indecomposable objects, and consequently so does $\mbox{mod}(Bf)$. Hence $fBf$ is representation-finite.

(2) For each $A$-module $Y$ in $\mbox{mod}(Ae)$, the assumption (2) guarantees that $Y$ is a direct summand of $FH(Y)$. Using the same arguments above, we can prove that $eAe$ is representation-finite provided that $fBf$ is representation-finite. $\square$

\medskip
 Note that Frobenius extensions with the above conditions (1) and (2) in Proposition \ref{PropalgExtension} appear frequently in stable equivalences of Morita type. In fact, by a result in \cite[Corollary 5.1]{DugasVilla}, if $A$ and $B$ are algebras such that their semisimple quotients are separable and if at least one of them is indecomposable, then there is a $k$-algebra $\Lambda$, Morita equivalent to $A$, and an injective ring homomorphism $B\hookrightarrow \Lambda$ such that $_{\Lambda}\Lambda\otimes_B\Lambda_{\Lambda}\simeq {}_{\Lambda}\Lambda_{\Lambda} \oplus {}_{\Lambda}P_{\Lambda}$ and $_B\Lambda_B\simeq {}_BB_B\oplus {}_BQ_B$ with $P$ and $Q$ are projective bimodules. This means that the extension $B\hookrightarrow \Lambda$ is a split, separable Frobenius extension.

\medskip
Let us mention a special case of Proposition \ref{PropalgExtension}. Suppose that $A$ is an algebra and $G$ is a finite group together with a group homomorphism from $G$ to $\mbox{Aut}(A)$, the group of automorphisms of the $k$-algebra $A$. Then one may form the skew group algebra $A*G$ of $A$ by $G$ and get the following corollary.

\begin{Koro}
Let $A$ be an algebra, and let $A*G$ be the skew group algebra of $A$ by $G$ with $G$ a finite group. If the order of $G$ is invertible in $A$, then $A*G$ is Frobenius-finite if and only if so is $A$.
\end{Koro}

{\it Proof.} Note that $A$ is a subalgebra of $A*G$. We just need to verify all the conditions in Proposition \ref{PropalgExtension}. However, all of them follow from \cite[Theorem 1.1]{ReitenRiedtmann}. $\square$

\medskip
\medskip
Next, we shall show that cluster-tilted algebras are Frobenius-finite. Suppose that $H$ is a finite-dimensional hereditary algebra over an algebraically closed field. Let $\tau_{D}$ be the Auslander-Reiten translation functor on $\Db{H}$, and let  ${\cal C}$ be  the orbit category $\Db{H}/\langle\tau_D^{-1}[1]\rangle$, which is a triangulated category with Auslander-Reiten translation $\tau_{\cal C}$.  Let ${\cal S}$ be the class of objects in $\Db{H}$ consisting of all modules in  $\modcat{H}$ and the objects $P[1]$, where $P$ runs over all modules in $\pmodcat{H}$.  The following two facts follows from \cite[Propositions 1.3, 1.6]{BMRRT}.

\medskip
  (a) $\tau_DX$ and $\tau_{\cal C}X$ are isomorphic in ${\cal C}$ for each object $X$ in $\Db{H}$;

  (b) Two objects $X$ and $Y$ in ${\cal S}$ are isomorphic in ${\cal C}$ if and only if they are isomorphic in $\Db{H}$;

 (c) $\Hom_{\cal C}(X, Y)=\Hom_{\Db{H}}(X, Y)\oplus \Hom_{\Db{H}}(X, \tau_D^{-1}Y[1])$ for all $X, Y\in {\cal S}$.  In particular, for each $H$-module $X$, then $\End_{\cal C}(X)=\End_{\Db{H}}(X)\ltimes\Hom_{\Db{H}}(X, \tau_D^{-1}X[1])$, the trivial extension of $\End_{\Db{H}}(X)$  by the bimodule $\Hom_{\Db{H}}(X, \tau_D^{-1}X[1])$ (see \cite[Proposition 1.5]{BMRRT}).

 Recall that, given an algebra $A$  and an $A$-$A$-bimodule $M$, the \emph{trivial extension} of $A$ by $M$, denoted by $A\ltimes M$, is the algebra with the underlying $k$-module $A\oplus M$ and the multiplication given by
$$(a, m)(a', m'):=(aa', am'+ma') \mbox{\; for }\; a, a'\in A, \, m, m'\in M.$$
If $M=DA$, then $A\ltimes DA$ is simply called the trivial extension of $A$, denoted by $\mathbb{T}(A)$.

\medskip
If $T$ is a cluster-tilting object in ${\cal C}$, then its endomorphism algebra $\End_{\cal C}(T)$ is called a \emph{cluster-tilted algebra}. Let $T$ be a basic tilting $H$-module. Then $\End_{\cal C}(T)$ is a cluster-tilted algebra and all cluster-tilted algebras can be obtained in this way.

Recall that the modules in $\add\{\tau_H^{-i}H|i\geq 0\}$ are called preprojective modules, and the modules in $\add\{\tau_H^i D(H)\mid i\geq 0\}$ are called preinjective modules.

\begin{Prop}
All cluster-tilted algebras  are Frobenius-finite.\label{cluster}
\end{Prop}

{\it Proof.} Let $A$ be a cluster-tilted algebra. Then, without loss of generality, we assume that $A=\End_{\cal C}(T)$, where $T$ is a basic tilting module over a connected, finite-dimensional hereditary algebra $H$ over an algebraically closed field.  If $H$ is of Dynkin type, then $A$ is  representation-finite and, of course, Frobenius-finite.

\smallskip
From now on, we assume that $H$ is representation-infinite.  Using a method similar to the one in the proof of \cite[Lemma 1]{RingelSelfInjCluster}, we deduce that the associated self-injective algebra of $A$ is isomorphic to $\End_{\cal C}(T')$  where $T'$ is a maximal direct summand of $T$ with $\tau_{\cal C}^2T'\simeq T'$ in ${\cal C}$.  By the fact (a) above, the objects $\tau_{D}^2T'$ and $T'$ are isomorphic in ${\cal C}$.  Suppose that $T'$ has a decomposition $T'=U\oplus M\oplus E$ such that $U$ is preprojective, $M$ is regular and $E$ is preinjective. For each projective $H$-module $P$, we have an Auslander-Reiten triangle
$$\nu_HP[-1]\lra V\lra P\lra \nu_HP$$
in $\Db{H}$, showing that $\tau_DP=\nu_HP[-1]$. Thus  $\tau_D^2P$, which is just $\tau_D(\nu_HP)[-1]$,   is isomorphic in ${\cal C}$ to $\nu_HP$ since ${\cal C}$ is the orbit category of $\Db{H}$ with respect to the auto-equivalence functor $\tau_D^{-1}[1]$.  Since $H$ is representation-infinite, for each $i>0$,  the object $\tau_D^i(\nu_HP)$ is isomorphic in $\Db{H}$ to $\tau_H^i(\nu_HP)$ which is a preinjective $H$-module.  Hence $\tau_D^mP$ is isomorphic in ${\cal C}$ to a preinjective $H$-module for all $m\geq 2$.  It follows that, for each preprojective $H$-module $V$,  the object $\tau_D^nV$ is isomorphic in ${\cal C}$ to a preinjective module provided $n$ is big enough.  Applying $\tau_D$ to a regular (preinjective, respectively) $H$-module always results in a regular (preinjective, respectively) $H$-module.  Thus, by applying $\tau_D^{2n}$ with $n$ large enough,  $\tau_D^{2n}T'\simeq \tau_D^{2n}U\oplus \tau_D^{2n}M\oplus \tau_D^{2n}E$ is isomorphic in ${\cal C}$ to an $H$-module $T''$ which has no preprojective direct summands.  Hence $T'$ and $T''$ are isomorphic in ${\cal C}$. By the fact (b) above, $T'$ and $T''$ are isomorphic in $\Db{H}$, and therefore they are also isomorphic as $H$-modules. However $T''$ has no preprojective direct summands. This forces $U$ to be zero.  Dually, one can prove that $E=0$.  Hence $T'$ is actually a regular $H$-module.  In this case $\tau_D^2T'$ is just $\tau_H^2T'$. By the fact (b) again, we see that $\tau_H^2T'$ and $T'$ are isomorphic in $\Db{H}$, and consequently $\tau_H^2T'\simeq T'$ as $H$-modules.

If $H$ is wild, then there is
no $\tau_H$-periodic $H$-modules at all. Hence $T'=0$ and $A$ is Frobenius-free in this case.  If $H$ is tame, then we have the following algebra isomorphisms
$$\begin{array}{rl}
\End_{\cal C}(T') &=\End_{\Db{H}}(T')\ltimes \Hom_{\Db{H}}(T', \tau_{D}^{-1}T'[1]) \quad \mbox{(by the fact (c) above)}\\
  &\simeq  \End_{H}(T')\ltimes \Ext_H^1(T', \tau_{H}^{-1}T')\\
  &\simeq  \End_{H}(T')\ltimes D\Hom_H(\tau_H^{-1}T', \tau_{H}T')  \quad \mbox{(by Auslander-Reiten formula)}\\
  &\simeq  \End_{H}(T')\ltimes D\Hom_H(T', \tau_{H}^{2}T') \\
  &\simeq  \End_{H}(T')\ltimes D\Hom_H(T', T')\\
  &= \End_H(T')\ltimes D\End_H(T')
\end{array}$$
is the trivial extension of $\End_H(T')$.  We claim that $\mathbb{T}(\End_{H}(T'))$ is
representation-finite.  Since $T$ is a tilting module over a tame
hereditary algebra $H$, it must contain either an indecomposable
preprojective or preinjective summand (see, for example,  the proof
of \cite[Lemma 3.1]{HappelRingel}). Thus there is an integer $n$
with $|n|$ minimal, such that $\tau_H^nT$  has a non-zero projective
or injective direct summand.  Assume that $\tau_H^nT\simeq He\oplus
X$  for some idempotent $e$ in $H$, and $X$ has no projective direct summands. Then $\tau_HX$ is
a tilting $H/HeH$-module.  Thus $\End_{H}(X)\simeq
\End_{H}(\tau_HX)$ is a tilted algebra of Dynkin type (not
necessarily connected), and consequently its trivial extension
$\mathbb{T}(\End_{H}(X))$ is representation-finite (see
\cite[Chapter V]{HappelTriangle}).  Since $T'$ is $\tau_H$-periodic,
$\tau_H^nT'$ has to be a direct summand of $X$. Thus,
$\End_H(T')\simeq \End_H(\tau_H^nT')$ is isomorphic to
$f\End_{H}(X)f$ for some idempotent $f$ in $\End_H(X)$. Hence the
trivial extension $\mathbb{T}(\End_H(T'))$ is isomorphic to
$f\mathbb{T}(\End_H(X))f$, and therefore it is
representation-finite.  When $\tau_H^nT$ contains an injective
direct summand, the proof can be proceeded similarly. $\square$

%\begin{Def} A $k$-algebra $A$ is said to be almost-Frobenius
%if the basic algebra of $A$ has at most one indecomposable projective
%module (up to isomorphism) which is not $\nu$-stable.\end{Def}

\subsection{Proof of Theorem \ref{ThmRepFin}}
Throughout this subsection, $k$ denotes an algebraically closed field.
The main idea of the proof of Theorem \ref{ThmRepFin} is
to utilize Theorem \ref{ThmLift} inductively. The following lemma is crucial to the induction procedure.

\begin{Lem} Let $A$ and $B$ be representation-finite,
self-injective $k$-algebras without semisimple direct summands.
Suppose that $\Phi: \stmodcat{A}\ra\stmodcat{B}$ is a stable
equivalence of Morita type. Then there is a simple $A$-module $X$
and integers $r$ and $t$ such that $\tau^r\circ \Omega^t\circ \Phi
(X)$ is isomorphic in $\stmodcat{B}$ to a simple $B$-module, where $\tau$ and $\Omega$ stands for the Auslander-Reiten translation and Heller operator, respectively.
\label{LemRepFinSimple}
\end{Lem}

{\it Proof.} Let $\Gamma_s(A)$ denote the stable Auslander-Reiten
quiver of $A$ which has isomorphism classes of non-projective
indecomposable $A$-modules as vertices and irreducible maps as
arrows. Then $\Gamma_s(A)$ and $\Gamma_s(B)$ are isomorphic as translation quivers. By \cite{LiuSummands}, we may assume that the algebras $A$ and $B$ are indecomposable. Then $\Gamma_s(A)$ and $\Gamma_s(B)$ are of the form
$\mathbb{Z}\Delta/G$ for some Dynkin graph $\Delta=A_n, D_n(n\geq
4), E_n(n=6,7,8)$ and a nontrivial admissible automorphism group $G$
of $\mathbb{Z}\Delta$ (\cite{Riedt}). We fix an isomorphism $s_A:
\mathbb{Z}\Delta/G\ra\Gamma_s(A)$, and set
$$\pi_A: \mathbb{Z}\Delta\lraf{\rm can}\mathbb{Z}\Delta/G\lraf{s_A}\Gamma_s(A).$$
Then $\pi_A$ is a covering map of translation quivers (see
\cite{Riedt}). Now we fix some isomorphisms of these translation
quivers.

\smallskip
$\bullet$ The Heller operator $\Omega_A$ gives rise to an automorphism $\omega_A: \Gamma_s(A)\ra \Gamma_s(A)$.

$\bullet$ The Auslander-Reiten translation $\tau_A$ gives rise to an automorphism $\tau_A: \Gamma_s(A)\ra \Gamma_s(A)$.

$\bullet$ Similarly, we have two automorphisms $\omega_B$ and $\tau_B: \Gamma_s(B)\ra \Gamma_s(B)$.

$\bullet$ The functor $\Phi$ induces an isomorphism $\phi: \Gamma_s(A)\ra \Gamma_s(B)$.

\smallskip
{\parindent=0pt Since the stable equivalence $\Phi$ is of Morita
type}, we have $\tau_A\phi=\phi\tau_B$ and
$\omega_A\phi=\phi\omega_B$. We set $\pi_B:=\pi_A\phi$. Then $\pi_B$
is also a covering map.

Let $\Delta$ be a Dynkin diagram of $n$ vertices. For the vertices of $\mathbb{Z}\Delta$, we use the coordinates $(s,t)$ with $1\leq t\leq n$ as described in \cite[fig. 1]{BrLaserRiedt}. A vertex $(p, 1)$ with $p\in\mathbb{Z}$ is called a {\em bottom vertex}. The vertices $(p, n)$ in $\mathbb{Z}A_n$ and $(p, 5)$ in $\mathbb{Z}E_6$ with $p\in\mathbb{Z}$ are called {\em top vertices}.

By definition, $\tau_{\Delta}: (p,q)\mapsto (p-1, q)$ is the translation on $\mathbb{Z}\Delta$ and all homomorphisms of translation quivers commute with the translation. The automorphism $\omega_A$ can be lifted to an admissible automorphism $\omega_{\Delta}$ of $\mathbb{Z}\Delta$ such that $\pi_A\omega_A=\omega_{\Delta}\pi_A$. The automorphism $\omega_{\Delta}$ can be defined as follows:  If $\Delta=A_n$, then $\omega_{A_n}(p,q)=(p+q-n, n+1-q)$ (see \cite[Section 4]{IyamaHighAR}). Using the method in \cite[4.4]{IyamaHighAR}, one can easily get that $\omega_{E_6}(p,q)=(p+q-6,6-q)$ for $q\neq 6$, and $\omega_{E_6}(p,6)=(p-6,6)$. Note that the method in \cite[4.4]{IyamaHighAR} does not depend on higher Auslander-Reiten theory and its main ingredients are actually the Auslander-Reiten formula and ordinary Auslander-Reiten theory.  Thus, for $\Delta=A_n$ or $E_6$, the automorphism $\omega_{\Delta}$ interchanges top vertices and bottom vertices.

Let ${\cal S}_A$ and ${\cal S}_B$ be complete sets of isomorphism classes of simple modules over $A$ and $B$, respectively. Define $\mathscr{C}_A:=\{x\in\mathbb{Z}\Delta|(x)\pi_A\in {\cal S}_A\}$ and $\mathscr{C}_B:=\{x\in\mathbb{Z}\Delta|(x)\pi_B\in {\cal S}_B\}$. Note that $\mathscr{C}_A$ and $\mathscr{C}_B$ are ``configurations" on $\mathbb{Z}\Delta$ by \cite[Propositions 2.3 and 2.4]{RiedtAn}). For the precise definition of configurations, we refer the reader to \cite{RiedtAn}. Note that if $\mathscr{C}$ is a configuration on $\mathbb{Z}\Delta$, then so is the image $(\mathscr{C})g$ for any admissible automorphism $g$ of $\mathbb{Z}\Delta$. In particular, $(\mathscr{C})\omega_{\Delta}$ and $(\mathscr{C})\tau_{\Delta}$ are configurations for all configurations $\mathscr{C}$.

\medskip
Claim 1: Each configuration $\mathscr{C}$ on $\mathbb{Z}{A_n}$ contains either a top vertex or a bottom vertex.

\medskip
{\it Proof.} Recall from \cite[Proposition 2.6]{RiedtAn} that there is a bijection between the configurations on $\mathbb{Z}A_n$ and the partitions $\sigma$ of the vertices of the regular $n$-polygon such that the convex hulls of different parts of $\sigma$ are disjoint. For such a partition $\sigma$, it is easy to see that either there is a part consisting of a single vertex, or there is a part containing two adjoint vertices. By using the bijection \cite[Proposition 2.6]{RiedtAn}, we can see that in the former case, the corresponding configuration contains a vertex $(i,n)$ for some integer $i$, and in the latter case, the corresponding configuration contains $(j,1)$ for some integer $j$. $\square$

\medskip
Claim 2: Let $\mathscr{C}$ be a configuration on $\mathbb{Z}\Delta$ with $\Delta= A_n, D_n(n\geq 4), E_6, E_7$ or $E_8$. Then either $\mathscr{C}$ or $(\mathscr{C})\omega_{\Delta}$ contains a bottom vertex.

\medskip
{\it Proof.}  We verify the statement in several cases.

(a) $\Delta=A_n$. Since $\omega_{A_n}$ maps top vertices to bottom vertices, Claim 2 follows from Claim 1.

(b) $\Delta=D_n$. The statement for $\mathbb{Z}D_4$ follows directly from \cite[7.6]{BrLaserRiedt}. Suppose $n\geq 5$. For $m\leq n-2$, let $\psi_m: \mathbb{Z}A_m\ra \mathbb{Z}D_n$ be the embedding defined in \cite[Section 6]{RiedtConfigDn}. By definition, $\psi_m$ maps all top and bottom vertices of $\mathbb{Z}A_m$ to bottom vertices of $\mathbb{Z}D_n$. By the two propositions in \cite[Section 6]{RiedtConfigDn}, each configuration on $\mathbb{Z}D_n$ contains the image of some configuration on $\mathbb{Z}A_r$ under $\tau_{D_n}^t\psi_r$ for some $0<r\leq n-2$ and $t\in\mathbb{Z}$. Together with Claim 1, this implies that each configuration on $\mathbb{Z}D_n$ with $n\geq 5$ contains at least one bottom vertex.

(c) $\Delta=E_6$. Note that $\omega_{E_6}$  maps top vertices to bottom vertices, and all the automorphisms of $\mathbb{Z}E_6$ are of the form $\tau_{\Delta}^s\omega_{\Delta}$ for some integer $s$ ( see \cite{Riedt}). Thus, the claim for $E_6$ follows from the list of isomorphism classes of configurations on $\mathbb{Z}E_6$ given in \cite[Section 8]{BrLaserRiedt}

(d) $\Delta=E_7$ or $E_8$. All the automorphisms of $\mathbb{Z}E_7$ and $\mathbb{Z}E_8$  are of the form $\tau_{\Delta}^s$ for some integer $s$. The claim then follows by checking the list of isomorphism classes of configurations on $\mathbb{Z}E_7$ and $\mathbb{Z}E_8$ given in \cite[Section 8]{BrLaserRiedt}.  $\square$

\medskip
Using Claim 2, we can assume that $(\mathscr{C}_A)\omega_{\Delta}^{a}$ contains a bottom vertex $(r_1, 1)$, and that $(\mathscr{C}_B)\omega_{\Delta}^{b}$ contains a bottom vertex $(r_2, 1)$, where $a,b$ are taken from $\{0, 1\}$. Let $x$ be
in $\mathscr{C}_A$ such that $(x)\omega_{\Delta}^a=(r_1,1)$. Then $(x)\omega_{\Delta}^a\tau_{\Delta}^{(r_1-r_2)}=(r_2,1)$, and
$$y:=(x)\omega_{\Delta}^{(a-b)}\tau_{\Delta}^{(r_1-r_2)}=(x)\omega_{\Delta}^{a}\tau_{\Delta}^{(r_1-r_2)}\omega_{\Delta}^{-b}=(r_2,1)\omega_{\Delta}^{-b}\in \mathscr{C}_B$$
Let $r=r_1-r_2$ and $t=a-b$. Then
$$(x)\pi_A\phi\omega_B^t\tau_B^r=(x)\pi_A\omega_A^t\tau_A^r\phi=(x)\omega_{\Delta}^t\tau_{\Delta}^r\pi_A\phi=(y)\pi_B.$$
Thus, the simple $A$-module $X:=(x)\pi_A$ is sent to the simple $B$-module $Y:=(y)\pi_B$ by the functor $\tau_B^{r}\circ\Omega_B^t\circ\Phi$
up to isomorphism in $\stmodcat{B}$.  $\square$

It would be nice to have a homological proof of Lemma \ref{LemRepFinSimple}.

\medskip
We have now accumulated all information necessary to prove the main result Theorem \ref{ThmRepFin}.

\medskip
{\bf Proof of  Theorem \ref{ThmRepFin}.}

Let $\Phi: \stmodcat{A}\ra \stmodcat{B}$ be a stable equivalence of Morita type. Suppose that $\Delta_A$ and $\Delta_B$ are the associated self-injective algebras of $A$ and $B$, respectively. Then it follows from \cite[Theorem 4.2]{DugasVilla} that $\Phi$ restricts to a stable equivalence $\Phi_1: \stmodcat{\Delta_A}\ra \stmodcat{\Delta_B}$ of Morita type. By Theorem \ref{ThmLift}, the stable equivalence $\Phi$ lifts  derived equivalence provided that $\Phi_1$ lifts to a derived equivalence. By the definition of associated self-injective algebras and Lemma \ref{LemIdemProp} (4), the algebras $\Delta_A$ and $\Delta_B$ have no semisimple direct summands.

If $\Delta_A=0$, then $\Phi$ lifts to a Morita equivalence between $A$ and $B$, and therefore Theorem \ref{ThmRepFin} follows. So we may suppose that $\Delta_A$ is not zero. Then, by Lemma \ref{LemRepFinSimple}, there are integers $r$ and $s$ such that the functor $\tau^r\Omega^s\Phi_1: \stmodcat{\Delta_A}\ra \stmodcat{\Delta_B}$ sends some simple $\Delta_A$-module to some simple $\Delta_B$-module. If the numbers of non-isomorphic simple modules over $\Delta_A$ and $\Delta_B$ equal $1$, then Proposition \ref{PropLiftMorita} provides a Morita equivalence between $\Delta_A$ and $\Delta_B$. So we may assume that $\Delta_A$ and $\Delta_B$ have more than $1$ simple modules. In this case, we can find $\nu$-stable idempotent elements $e$ and $f$ in $\Delta_A$ and $\Delta_B$, respectively, such that the equivalence $\tau^r\Omega^s\Phi_1$ restricts to a stable equivalence $\Phi_2: \stmodcat{e\Delta_Ae}\ra \stmodcat{f\Delta_Bf}$ of Morita type, and that the algebras $e\Delta_Ae$ and $f\Delta_Bf$  have less mumber of non-isomorphic simple modules than $\Delta_A$ and $\Delta_B$ do, respectively. Since $e\Delta_Ae$ and $f\Delta_Bf$ are again representation-finite and self-injective without semisimple direct summands, we can assume, by induction, that $\Phi_2$ lifts to a derived equivalence. Thus, by Theorem \ref{ThmLift}, the stable equivalence $\tau^r\Omega^s\Phi_1$ lifts to a derived equivalence. Moreover, for self-injective algebras, both $\tau$ and $\Omega$ lift to derived equivalences between $\Delta_A$ and $\Delta_B$. Hence $\Phi_1$ lifts to a derived equivalence. $\square$

\medskip
{\it Remark}. For standard representation-finite, self-injective $k$-algebras $A$ and $B$ not of type $(D_{3m}, s/3,1)$ with $m\geq 2$ and $3\nmid s$, Asashiba proved in \cite{AsashibaLift} that every individual stable equivalence between $A$ and $B$ lifts to a derived equivalence. This was done by his derived equivalence classification of representation-finite, self-injective algebras.  In Theorem \ref{ThmRepFin} we consider instead stable equivalences of Morita type, and in this case, we can deal with all representation-finite, self-injective algebras without care about the types. Also, the proof of Theorem \ref{ThmRepFin} is independent of Asashiba's derived equivalence classification of representation-finite, self-injective algebras. So we have the following generalization of Asashiba's result.

\begin{Koro} If $A$ and $B$ are arbitrary representation-finite self-injective algebras over an algebraically closed field without semisimple direct summands, then every stable equivalence of Morita type between $A$ and $B$ can be lifted to an iterated almost $v$-stable derived equivalence.\label{extensionofasshiba}
\end{Koro}

As another consequence of Theorem \ref{ThmRepFin}, we have

\begin{Koro} If $A$ and $B$ are the Auslander $k$-algebras without semisimple direct summands, then every individual stable equivalence of Morita type between $A$ and $B$ lifts to an iterated almost $\nu$-stable derived equivalence. \label{cor5.5}
\end{Koro}

{\it Proof.} By a result of Auslander (see, for example, \cite[Theorem 5.7]{AusReiten}), we may assume that $A$ is the endomorphism algebra of a representation-finite algebra $A'$. Thus
the Frobenius parts of $A$ has to be of the form $eA'e$ with $e^2=e\in A'$. Therefore it is representation-finite since so is $A'$. Thus $A$ is Frobenius-finite. Now Corollary \ref{cor5.5} follows immediately from Theorem \ref{ThmRepFin}. $\square$

%\section{Some examples of derived equivalent blocks of finite groups\label{sect6}}
\section{A machinery for lifting stable equivalences to derived equivalences\label{sect6}}

In this section, we give a procedure for lifting a class of stable equivalences of Morita type to derived equivalences. With this machinery we re-check some derived equivalent block algebras of finite groups.

Let $A$ be an algebra, and let ${\cal S}_A$ be a
complete set of pairwise non-isomorphic simple $A$-modules. For each
simple $A$-module $V\in {\cal S}_A$, we fix a primitive idempotent element
$e_V$ in $A$ with $e_V\cdot V\neq 0$, such that the idempotent elements $\{e_{V}\mid
V\in{\cal S}_A\}$ are pairwise orthogonal. Thus, for
any nonempty subset $\sigma$ of ${\cal S}_A$, the element
$e_{\sigma}:=\sum_{V\in\sigma}e_V$ is an idempotent element in $A$.

\medskip
Theorem \ref{ThmLift} and the proof of Theorem \ref{ThmRepFin}
suggest an inductive method to check whether a stable equivalence of
Morita type can be induced by a derived equivalence. The procedure reads
as follows:

\medskip
{\bf Assumption}: Let $\Phi:\stmodcat{A}\ra \stmodcat{B}$ be a
stable equivalence of Morita type between two algebras without semisimple direct summands.
%\textcolor[rgb]{0.98,0.00,0.00}
Suppose that $A/\rad(A)$ and $B/\rad(B)$ are separable.

\smallskip
{\bf Step $1$:} If There is a simple $A$-module $V$ such that $\Phi(V)$ is a simple $B$-module, then we set
$$\sigma:=\{V\in {\cal S}_A \mid \Phi(V) \mbox{
is non-simple}\}\mbox{ and } \sigma':={\cal S}_B\backslash
\Phi({\cal S}_A\backslash \sigma).$$ By Lemma \ref{LemRes}, the
functor $\Phi$ restricts to a stable equivalence of Morita type
between $e_{\sigma}Ae_{\sigma}$ and $e_{\sigma'}Be_{\sigma'}$. Moreover, the idempotent elements $e_{\sigma}$ and $e_{\sigma'}$ are both $\nu$-stable. In fact, by Lemma \ref{LemPropStM} (5) and (6),
for each $V$ in ${\cal S}_A$, the $B$-module $\Phi(V)$ is non-simple if and only if $\Hom_A({}_AP, V)\neq 0$, or equivalently, $V\in\add(\top({}_AP))$, where $P$ is given in the definition of the stable equivalence $\Phi$ of Morita type.  This implies that $\add(Ae_{\sigma})=\add({}_AP)$.  It follows from Lemma \ref{LemPropStM} (2) that $e_{\sigma}$ is $\nu_A$-stable. Similarly, it can be shown that $e_{\sigma'}$ is $\nu_B$-stable. By Lemma \ref{LemIdemProp} (3), the algebras $e_{\sigma}Ae_{\sigma}$ and
$e_{\sigma'}Be_{\sigma'}$ are self-injective with less simple
modules.

\smallskip
{\bf Step $2$:} Find some suitable stable equivalence $\Xi: B\ra C$
of Morita type between the algebra $B$ and another algebra $C$, which is induced by a derived equivalence such that the composite
$\Xi\circ\Phi$ sends some simple $A$-modules to simple $C$-modules.
Then go back to Step $1$. Once we get two representation-finite algebras in the procedure, Theorem \ref{ThmRepFin} will be applied.

\medskip
This procedure is somewhat similar to, but different from the method
of Okuyama in \cite{Oku1997}: in our procedure, Step
$1$ always reduces the number of simple modules and makes
the situation considered easier after each step, while the procedure in \cite{Oku1997} does not change the number of simple modules.

\medskip
In the following, we will illustrate the
ideas mentioned above by examples.

\medskip
{\bf\parindent=0pt Example 1:} In \cite{JMM}, it was proved that the
Brou\'{e}'s Abelian Defect Group Conjecture is true for the faithful $3$-blocks of defect $2$
of $4.M_{22}$, which is the non-split central extension of the
sporadic simple group $M_{22}$ by a cyclic group of order 4.
Now we shall show that the procedure described above can be used to give
a short proof of the conjecture in this case, which avoids many technical calculations, comparing with the original proof in \cite{JMM}.

It is known that each of the two block algebras $B_{+}$ and $b_{+}$ has 5 simple
modules. The simple $B_{+}$-modules are labeled by $56a, 56b, 64, 160a,
160b$, and the simple $b_{+}$-modules are labeled by $1a, 1b, 2, 1c$
and $1d$. There is a stable equivalence
$$\Phi: \stmodcat{B_{+}}\lra \stmodcat{b_{+}}$$ of Morita type (see \cite{JMM}) such that
$$\Phi(56a)=\Omega^{-1}(1a), \Phi(56b)=\Omega(1b), \Phi(160a)=1c, \Phi(160b)=1d,$$
and $\Phi(64)$ has the following Loewy structure
$$\left[\begin{matrix}1b\\2\\1a\end{matrix}\right].$$
For $x\in\{a,b,c,d\}$ and $\{y, y',
y''\}=\{a,b,c,d\}\backslash \{x\}$, the Loewy structures of the
projective $b_{+}$-modules are
$$P(1x): \left[\begin{matrix}1x\\2\\1y\,\, 1y'\,\,
1y''\\2\\1x\end{matrix}\right],\quad  P(2):
\left[\begin{matrix}2\\1a \,\, 1b \,\, 1c \,\, 1d\\ 2 \, \,
2 \,\, 2\\ 1a \,\, 1b \,\, 1c \,\, 1d\\2\end{matrix}\right].$$
Now, we use Steps $1$ and  $2$ repeatedly and verify
that the stable equivalence $\Phi$ lifts to a derived equivalence.

Note that $\Phi$ sends the simple module $160b$ to a simple module. So we can use Step 1. Let $\sigma=\{56a, 56b, 64\}$, and $\sigma'=\{1a,1b,2\}$. Then
$\Phi$ restricts to a stable equivalence of Morita type
$$\Phi_1:
\stmodcat{e_{\sigma}B_+e_{\sigma}}\lra\stmodcat{e_{\sigma'}b_+e_{\sigma'}}.$$
The Loewy structures of the projective $e_{\sigma'}b_+e_{\sigma'}$-modules
$e_{\sigma'}P(1a)$ and $e_{\sigma'}P(1b)$ are
$$e_{\sigma'}P(1a): \left[\begin{matrix}1a\\2\\1b\\2\\1a\end{matrix}\right],\mbox{ \quad and \quad }
e_{\sigma'}P(1b):\left[\begin{matrix}1b\\2\\1a\\2\\1b\end{matrix}\right].$$
The images of the simple modules under $\Phi_1$  are
$$\Phi_1(56a)\simeq \left[\begin{matrix}1a\\2\\1b\\2\end{matrix}\right],
\Phi_1(56b)\simeq
\left[\begin{matrix}2\\1a\\2\\1b\end{matrix}\right], \mbox{ \quad and \quad
}\Phi_1(64)\simeq
\left[\begin{matrix}1b\\2\\1a\end{matrix}\right].$$ By
\cite{Oku1997}, the idempotent $e=e_{1a}+e_{1b}$ defines a tilting
complex $\cpx{T}$ over $e_{\sigma'}Ae_{\sigma'}$. Setting
$C:=\End(\cpx{T})$ and labeling the simple $C$-modules by $1a, 1b$
and $2$, the derived equivalence between $e_{\sigma'}Ae_{\sigma'}$
and $C$ induces a stable equivalence of Morita type $\Xi:
\stmodcat{e_{\sigma'}Ae_{\sigma'}}\ra \stmodcat{C}$ such that
$\Xi(2)\simeq 2$,
$\Xi(\left[\begin{smallmatrix}1b\\2\\1a\end{smallmatrix}\right])\simeq
1b$, and $\Xi(\left[\begin{smallmatrix}1a\\2\\1b\end{smallmatrix}\right])\simeq
1a$. Thus $\Xi\Phi_1(64a)\simeq 1b$, $\Xi\Phi_1(56a)\simeq
\genfrac{[}{]}{0pt}{}{1a}{2}$ and $\Xi\Phi_1(56b)\simeq
\genfrac{[}{]}{0pt}{}{2}{1a}$. Let $\sigma_1:=\{56a, 56b\}$ and
$\sigma'_1:=\{1a,2\}$. Then the composite $\Xi\Phi_1$ restricts to a
stable equivalence of Morita type
$$\Phi_2:\stmodcat{e_{\sigma_1}B_+e_{\sigma_1}}\lra
\stmodcat{e_{\sigma'_1}Ce_{\sigma'_1}}$$ such that
$\Phi_2(56a)=\left[\begin{smallmatrix}1a\\2\end{smallmatrix}\right]$
and
$\Phi_2(56b)=\left[\begin{smallmatrix}2\\1a\end{smallmatrix}\right]$.
Note that the Cartan matrix of $e_{\sigma'_1}Ce_{\sigma'_1}$ is
$\left[\begin{smallmatrix}2 &1\\1&3\end{smallmatrix}\right]$. It is easy to check that a symmetric algebra with this Cartan matrix is always representation-finite. Thus $\Phi_2$ lifts to a derived equivalence by Theorem \ref{ThmRepFin}, and consequently $\Phi$ lifts to a derived equivalence by our inductive procedure.
%However, we can do one step more as follows.
%The first column of the above Cartan matrix is the dimension vector of $e_{\sigma'_1}Ce_{1a}$, and indicates that the indecomposable projective module $e_{\sigma'_1}Ce_{1a}$ is of the form $\left[\begin{smallmatrix}1a\\2\\1a\end{smallmatrix}\right]$.
%This
%implies that $\Phi_2(56a)=\Omega^{-1}(1a)$, and
%$\Omega\Phi_2$ restricts to a stable equivalence $\Phi_3$ of Morita
%type between $e_{56b}B_+e_{56b}$ and $e_2Ce_2$, which sends $56b$ to $\genfrac{[}{]}{0pt}{}{2}{2}$. The algebra
%$e_2Ce_2$ is a local symmetric algebra of dimension $3$. Hence
%$\Omega\Phi_3$ sends $56b$ to $2$ and lifts to a
%Morita equivalence. Since $\Omega$ lifts to a derived equivalence,
%we deduce that $\Phi_3$ lifts to a derived equivalence. By our inductive method, the stable equivalence $\Phi$ lifts to
%a derived equivalence. 
The whole procedure can be illustrated by the following commutative diagram
$$\xymatrix@R=7mm@C=12mm{
\stmodcat{B_+} \ar[r]^{\Phi} & \stmodcat{b_+}\\
\stmodcat{e_{\sigma}B_+e_{\sigma}} \ar[r]^{\Phi_1}\ar[u]^{\lambda} &
\stmodcat{e_{\sigma'}b_+e_{\sigma'}}\ar[u]^{\lambda}\ar[r]^{\Xi} &
\stmodcat{C}\\
\stmodcat{e_{\sigma_1}B_+e_{\sigma_1}}\ar[u]^{\lambda}\ar[rr]^{\Phi_2}
    && \stmodcat{e_{\sigma'_1}Ce_{\sigma'_1}}\ar[u]^{\lambda} \\
    }$$
 with $\Phi_3$ lifting to a derived equivalence.

\medskip
{\parindent=0pt\bf Example 2:} Let $G$ be the Harada-Norton simple
group ${\bf HN}$, and let $k$ be an algebraically closed field of
characteristic $3$. In \cite{HNgroup},  the Brou\'{e}'s Abelian Defect Group Conjecture was
verified for non-principal blocks  of $kG$ with defect group
$C_3\times C_3$. In the following,
we will show how our results can be applied to give another proof to the conjecture in this case.
In fact, the two block algebras $A$ and $B$ have $7$ non-isomorphic
simple modules with ${\cal S}_A=\{1,2,3,4,5,6,7\}$ and ${\cal
S}_B=\{9a, 9b, 9c, 9d, 18a, 18b, 18c\}$, and there is a stable
equivalence $F: \stmodcat{A}\ra \stmodcat{B}$ of Morita type such
that
$$F(1)\simeq 9a, \quad F(2)\simeq 9b, \quad F(3)\simeq 9c$$
$$F(4)\simeq \xy (0,7)*+{18a}="u", (0, -7)*+{18a}="d",(7,0)*+{18c}="r",
(-7,0)*+{18b}="l",{\ar@{-} "u";"l"},{\ar@{-} "u";"r"},{\ar@{-} "d";"l"},
{\ar@{-} "d";"r"},\endxy, \quad F(5)\simeq \xy (0,7)*+{18c}="u", (0, -7)*+{18b}="d",
(7,0)*+{9d}="r", (-7,0)*+{9a}="l",{\ar@{-} "u";"l"},{\ar@{-} "u";"r"},{\ar@{-} "d";"l"},
{\ar@{-} "d";"r"},\endxy, \quad F(6)\simeq \xy (0,7)*+{18a}="u", (0, -7)*+{18a}="d",
(-7,0)*+{18c}="l", (7,0)*+{18b}="r",(14,7)*+{9d}="9du",(-14,-7)*+{9d}="9db",
{\ar@{-} "u";"l"},{\ar@{-} "u";"r"},{\ar@{-} "d";"l"},{\ar@{-} "d";"r"},
{\ar@{-} "r"; "9du"}, {\ar@{-} "l"; "9db"},\endxy, \quad F(7)
\simeq \xy (0,7)*+{18b}="u", (0, -7)*+{18c}="d",(7,0)*+{9c}="r",
(-7,0)*+{9b}="l",{\ar@{-} "u";"l"},{\ar@{-} "u";"r"},{\ar@{-} "d";"l"},{\ar@{-} "d";"r"},\endxy.$$
The Loewy structures of the indecomposable projective $B$-modules
$P(9d), P(18a), P(18b)$ and  $P(18c)$ are as follows.
$$P(9d): \left[\begin{matrix}9d\\18b\\9c\,\,18a\\18c\\9d\end{matrix}\right],\quad P(18a):
\xy
<16pt,0pt>:
   @={(0,0),(1,0),(2,0),(3,0),(4,0),(1,1),(3,1),(2,2),(1,-1),(3,-1),(2,-2)},
  s0*+{18a}="soc",
  s1*+{18b}="dr",
  s2*+{18c}="dl",
  s3*+{18a}="top",
  s4*+{18c}="ur",
  s5*+{18b}="ul",
  s6*+{9d}="mrr",
  s7*+{9a}="mr",
  s8*+{18a}="mm",
  s9*+{9c}="ml",
  s{10}*+{9b}="mll",
  {\ar@{-} "soc";"dl"},
  {\ar@{-} "soc";"dr"},
  {\ar@{-} "dr";"mrr"},
  {\ar@{-} "dr";"mr"},
  {\ar@{-} "dr";"mm"},
  {\ar@{-} "dl";"mll"},
  {\ar@{-} "dl";"ml"},
  {\ar@{-} "dl";"mm"},
  {\ar@{-} "ul";"mll"},
  {\ar@{-} "ul";"ml"},
  {\ar@{-} "ul";"mm"},
  {\ar@{-} "ur";"mrr"},
  {\ar@{-} "ur";"mr"},
  {\ar@{-} "ur";"mm"},
  {\ar@{-} "top";"ul"},
  {\ar@{-} "top";"ur"},
\endxy, \quad P(18b): \left[\begin{matrix}18b\\9b\,\, 18a\,\,9c\\18c\,\,18b\,\,18c\\9a\,\,18a\,\,9d\\18b
\end{matrix}\right],\quad P(18c): \left[\begin{matrix}18c\\9a\,\, 18a\,\,9d\\18b\,\,18c\,\,18b\\9b\,\,18a\,\,9c\\18c
\end{matrix}\right]$$
Taking $\sigma=\{4,5,6,7\}$ and $\sigma'=\{9d, 18a, 18b, 18c\}$, we see from Step $1$ that the functor $F$ restricts to a stable
equivalence of Morita type
$$F_1: \stmodcat{e_{\sigma}Ae_{\sigma}}\lra\stmodcat{e_{\sigma'}Be_{\sigma'}}$$
such that
$$F_1(4)\simeq \xy (0,7)*+{18a}="u", (0,
-7)*+{18a}="d",(7,0)*+{18c}="r", (-7,0)*+{18b}="l",{\ar@{-}
"u";"l"},{\ar@{-} "u";"r"},{\ar@{-} "d";"l"},{\ar@{-}
"d";"r"},\endxy, \quad F_1(5)\simeq \left[\begin{matrix}
18c\\9d\\18b\end{matrix}\right], \quad F_1(6)\simeq \xy (0,7)*+{18a}="u",
(0, -7)*+{18a}="d",(-7,0)*+{18c}="l",
(7,0)*+{18b}="r",(14,7)*+{9d}="9du",(-14,-7)*+{9d}="9db", {\ar@{-}
"u";"l"},{\ar@{-} "u";"r"},{\ar@{-} "d";"l"},{\ar@{-}
"d";"r"},{\ar@{-} "r"; "9du"}, {\ar@{-} "l"; "9db"},\endxy,\quad
F_1(7)\simeq \left[\begin{matrix} 18b\\18c\end{matrix}\right].$$ The
idempotent element $e_{18a}$ in $B$ defines a tilting complex $\cpx{T}$ over
$e_{\sigma'}Be_{\sigma'}$ (\cite{Oku1997}). Set
$C:=\End(\cpx{T})$ and label the simple $C$-modules by $9d, 18a,
18b$ and $18c$. Then the derived equivalence between
$e_{\sigma'}Be_{\sigma'}$ and $C$ induces a stable equivalence of
Morita type $\Xi:
\stmodcat{e_{\sigma'}Be_{\sigma'}}\ra\stmodcat{C}$ such that
$\Xi(9d)\simeq 9d$, $\Xi(18b)\simeq 18b$, $\Xi(18c)\simeq 18c$, and
$\Xi F_1(4)\simeq 18a$. Taking $\sigma_1=\{5,6,7\}$ and
$\sigma'_1=\{9d,18b,18c\}$, the functor $\Xi F_1$ restricts to a stable
equivalence of Morita type
$$F_2: \stmodcat{e_{\sigma_1}Ae_{\sigma_1}}\lra\stmodcat{e_{\sigma'_1}Ce_{\sigma'_1}}$$
such that $F_2(5)\simeq \left[\begin{smallmatrix}
18c\\9d\\18b\end{smallmatrix}\right], F_2(6)\simeq
\left[\begin{smallmatrix} 9d\\9d\end{smallmatrix}\right]$ and
$F_2(7)\simeq \left[\begin{smallmatrix}
18b\\18c\end{smallmatrix}\right]$. Note that the Cartan matrix of
$e_{\sigma'_1}Ce_{\sigma'_1}$ is
$\left[\begin{smallmatrix}2&1&1\\1&2&1\\1&1&3\end{smallmatrix}\right]$,
where the columns are dimension vectors of the projective modules $e_{\sigma'_1}Ce_{18b}$,
$e_{\sigma'_1}Ce_{18c}$ and $e_{\sigma'_1}Ce_{9d}$, respectively.
Then $F_2(5)\simeq \Omega^{-1}(18c)$. Thus, taking
$\sigma_2=\{6,7\}$ and $\sigma'_2=\{18b,9d\}$, the functor $\Omega
F_2$ restricts to a stable equivalence of Morita type
$$F_3: \stmodcat{e_{\sigma_2}Ae_{\sigma_2}}\lra\stmodcat{e_{\sigma'_2}Ce_{\sigma'_2}}.$$
The Cartan matrix of $e_{\sigma'_2}Ce_{\sigma'_2}$ is $\left[\begin{smallmatrix}2 &1\\1&3\end{smallmatrix}\right]$. This implies that 
$e_{\sigma'_2}Ce_{\sigma'_2}$ is representation-finite and that $F_3$ lifts to a derived equivalence by Theorem \ref{ThmRepFin}. 
Hence $F$ lifts to a derived equivalence.

\medskip
Finally, we point out that our methods also work for the most examples given in \cite{Oku1997}.

\smallskip
Let us end this section by mentioning the following questions suggested by our main results.

{\bf Question 1.} Given a stable equivalence $\Phi$ of Morita type between two self-injective algebras such that $\Phi$ does not send any  simple modules to simple modules, under which conditions can $\Phi$ be lifted to a derived equivalence?

{\bf Question 2.} Find more other sufficient conditions for stable equivalences of Morita type between general finite-dimensional algebras to be lifted to derived equivalences.

{\bf Question 3.} Find more classes of algebras that are Frobenius-finite. For example, when is a cellular algebra Frobenius-finite.

\medskip
{\bf Acknowledgement.} The research works of both authors are partially supported by NSFC. The author C.C.X.
thanks  BNSF(1132005, KZ201410028033) for partial support, while the author W.H. is grateful to the Fundamental Research Funds for the Central Universities for partial support.

\bigskip
{\footnotesize
\bigskip Wei Hu, School of Mathematical Sciences, Beijing Normal
University, 100875 Beijing, People's Republic of  China, and
Beijing Center for Mathematics and Information Interdisciplinary Sciences, 100048 Beijing, China

{\tt Email: huwei@bnu.edu.cn}

\bigskip
Changchang Xi, School of Mathematical Sciences, BCMIIS, Capital Normal
University, 100048 Beijing, People's Republic of  China

{\tt Email: xicc@cnu.edu.cn}}

\bigskip
First version: March 20, 2013; Revised: May 21, 2014. Final version: December 20, 2014.

\end{document}